# Analysis and Control of a Non-Local PDE Traffic Flow Model

**Iasson Karafyllis[*], Dionysis Theodosis[**] and Markos Papageorgiou[**]**

[*]Dept. of Mathematics, National Technical University of Athens,
Zografou Campus, 15780, Athens, Greece,
emails: iasonkar@central.ntua.gr , iasonkaraf@gmail.com

[**] Dynamic Systems and Simulation Laboratory,
Technical University of Crete, Chania, 73100, Greece
(emails: dtheodosis@dssl.tuc.gr , markos@dssl.tuc.gr)

## Abstract

The paper provides conditions that guarantee existence and uniqueness of classical solutions for a non-local conservation law on a ring-road with possible nudging (or "look behind") terms. The obtained conditions are novel, as they are not covered by existing results in the literature. The paper also provides results which indicate that nudging can increase the flow in a ring-road and, if properly designed, can have a strong stabilizing effect on traffic flow. More specifically, the paper gives results which guarantee local exponential stability of the uniform equilibrium profile in the $L^2$ state norm even for cases where the uniform equilibrium profile in a ring-road without nudging is not asymptotically stable and the model admits density waves. The efficiency of the use of nudging terms is demonstrated by means of a numerical example.

## 1. Introduction

Non-local traffic flow models with Partial Differential Equations (PDEs) are based on extensions of the well-known Lighthill-Whitham-Richards model (LWR model, see [19, 22]), where the speed is given by a non-local term. These models fall into the class of non-local conservation laws (see [5]) and possess some different features compared to the LWR model. Arrhenius "look-ahead" terms were considered in [24, 17, 18] as a result of stochastic microscopic dynamics, and it was shown that such models can develop shocks (and shock waves) in finite time. On the other hand, the fact that human drivers and automated vehicles adjust the vehicle speed based on a perception of downstream density, rather than the local density, motivated some researchers to express the perceived density by means of non-local (convolution) terms. Such models were studied in [2, 3, 4, 8, 14, 15], and it was shown that they may be producing smooth solutions.

In the era of automated vehicles, the real-time information fed to each vehicle on a road is exploited for the appropriate adjustment of the speed of the vehicle. In contrast to manual driving, this information may include upstream density data, in addition to downstream density data. Note that human drivers base their driving decisions only on the perceived *downstream* traffic state, something that leads to the celebrated anisotropy principle in traffic flow modeling [7]. The possible beneficiary role of the use of upstream density data was pointed out in [21], where the effect of the upstream density data on the speed adjustment was termed as "nudging". Such an effect was also studied in [17] (without reference to automated vehicles), where an Arrhenius "look-behind" non-local term was used for the mathematical expression of the use of upstream density data.

The selection of the nudging term for automated vehicles can be considered as a feedback design problem. Data are fed into the automated vehicles, based on which the vehicles adjust their speed. In other words, the density profile of the road changes over time, and this change is fed back to each automated vehicle. From a mathematically perspective, the use of upstream density data should not be performed in an arbitrary way, but so as to satisfy conditions for existence and uniqueness of solutions, together with further requirements for the closed-loop system (e.g., stability, optimality, etc.). It should be noticed here that the feedback design problem for the expression of "nudging" or "look-behind" effect can be considered as a special feedback design problem for non-local, hyperbolic PDEs (see [6, 13, 16]). However, this specific feedback design problem is different from other traffic control problems studied in the literature (see [11, 12, 25, 26, 27]).

The present paper answers these questions for a ring-road. We first present conditions which guarantee existence and uniqueness of classical solutions for a non-local conservation law with possible nudging terms (Theorem 2.3 and Theorem 2.4). The obtained results are novel, as they are not covered by the results in [2, 3, 4, 8, 14, 15], where either the use of upstream density data is not allowed or a ring-road is not studied. In addition, the present paper studies the effects of nudging and it is shown that:

(i) nudging can increase the flow in a ring-road at any density value;
(ii) if properly designed, nudging can have a strong stabilizing effect on ring-road traffic.

Indeed, we present results (Theorem 3.2) which guarantee local exponential stability of the uniform equilibrium profile in the $L^2$ state norm even for cases where the uniform equilibrium profile in a ring-road without nudging is not asymptotically stable and the model admits density waves. The existence of travelling waves for non-local conservation laws was studied in [23], where it is shown that travelling waves may occur even in non-local conservation laws.

The structure of the paper is as follows. Section 2 of the paper is devoted to the presentation of the non-local traffic flow models which are studied in the paper; moreover, the statements of the existence and uniqueness results for non-local traffic flow models are also given in Section 2. The effects of nudging on ring-road traffic are studied in Section 3. Illustrative numerical experiments are presented in Section 4, where the strong stabilizing effect of nudging is demonstrated. All proofs of the main results are provided in Section 5. Finally, concluding remarks are given in Section 6, and the Appendix contains the proof of an auxiliary result.

**Notation.** Throughout this paper, we adopt the following notation.

* $\Re_+ := [0,+\infty)$. For a real number $x \in \Re$, $[x]$ denotes the integer part of $x$, i.e., the greatest integer which is less or equal to $x$.
* For a vector $y \in \Re^N$, $y^\top$ denotes its transpose and $|y|_\infty = \max_{i=1,\dots,N}(|y_i|)$ denotes its infinity norm. For two vectors $x, y \in \Re^N$ we write $x \le y$ if and only if $x_i \le y_i$ for $i = 1,\dots,N$. The vector $1_N \in \Re^N$ is the vector $1_N = (1,\dots,1)^\top \in \Re^N$. We also define $y^{(0)} = y = (y_1,\dots,y_N)^\top \in \Re^N$, $y^{(1)} = (y_2,\dots,y_N,y_1)^\top \in \Re^N$ and $y^{(k)}$ for $k \ge 2$ by means of the recursive formula $y^{(k)} = (y^{(k-1)})^{(1)}$.
* Let $S \subseteq \Re^n$ be an open set and let $A \subseteq \Re^n$ be a set that satisfies $S \subseteq A \subseteq cl(S)$. By $C^0(A;\Omega)$, we denote the class of continuous functions on $A$, which take values in $\Omega \subseteq \Re^m$. By $C^k(A;\Omega)$, where $k \ge 1$ is an integer, we denote the class of functions on $A \subseteq \Re^n$, which takes values in $\Omega \subseteq \Re^m$ and has continuous derivatives of order $k$. In other words, the functions of class $C^k(A;\Omega)$ are the functions which have continuous derivatives of order $k$ in $S = \mathrm{int}(A)$ that can be continued continuously to all points in $\partial S \cap A$. When $\Omega = \Re$ then we write $C^0(A)$ or $C^k(A)$.
* Let $T \in (0,+\infty)$ and $\rho : [0,T] \times I \to \Re$ be given, where $I \subseteq \Re$ is an interval. We use the notation $\rho[t]$ to denote the profile at certain $t \in [0,T]$, i.e., $(\rho[t])(x) = \rho(t,x)$ for all $x \in I$. $L^p(I)$ with $p \ge 1$

denotes the equivalence class of measurable functions $f : I \to \Re$ for which $\|f\|_p = \left( \int_I |f(x)|^p dx \right)^{1/p} < +\infty$. $L^\infty(I)$ denotes the equivalence class of measurable functions $f : I \to \Re$ for which $\|f\|_\infty = \underset{x \in I}{ess\sup}\left(|f(x)|\right) < +\infty$. We use the notation $f'(x)$ for the derivative at $x \in I$ of a differentiable function $f : I \to \Re$.

* $W^{2,\infty}([0,1])$ is the Sobolev space of $C^1$ functions on $[0,1]$ with Lipschitz derivative.
* $Per(\Re)$ denotes the set of continuous, positive mappings $\rho : \Re \to (0,+\infty)$ which are periodic with period 1, i.e., $\rho(x+1) = \rho(x)$ for all $x \in \Re$.
* Let $(X, d_X)$ be a compact metric space and let $(Y, d_Y)$ be a given metric space. By $C^0(X\,;Y)$ we denote the set of continuous mappings $f : X \to Y$.

## 2. Non-Local Traffic Flow Models

Many non-local PDE traffic flow models which have appeared in the literature (see [2, 3, 4, 8, 14, 15]) have the form

$$\frac{\partial \rho}{\partial t}(t,x) + \frac{\partial}{\partial x}\left(\rho(t,x)v(t,x)\right) = 0, \text{ for } t \geq 0, x \in \Re \tag{2.1}$$

$$v(t,x) = f\left( \int_x^{x+\eta} \omega(s-x)\rho(t,s)ds \right), \text{ for } t \geq 0, x \in \Re \tag{2.2}$$

where $\rho(t,x)$ denotes the traffic density, $v(t,x)$ denotes the mean speed, $t \geq 0$ is time, $x$ is the spatial variable, $\eta > 0$ is a constant (reflecting the visibility area), $f : \Re_+ \to \Re_+$ and $\omega : \Re_+ \to \Re_+$ are non-increasing functions with $\int_0^\eta \omega(x)dx = 1$. Model (2.1), (2.2) constitutes a generalization of the classical LWR traffic flow model, where $f : \Re_+ \to \Re_+$ is the function that relates density to speed (fundamental diagram) and $\int_x^{x+\eta} \omega(s-x)\rho(t,s)ds$ is the downstream density perceived by the human driver at spatial position $x$. Thus, the driver adapts the speed according to (2.2) on the basis of the perceived downstream density.

As a farther generalization, when automated vehicles are present on a highway, there may be a benefit by allowing the vehicle speed to depend on upstream density levels as well. Such an effect has been termed in the literature as "nudging" (see [21]) or "look-behind" effect (see [17]). In this case, the speed may be given by a relation of the form

$$v(t,x) = f\left( \int_x^{x+\eta} \omega(s-x)\rho(t,s)ds \right) g\left( \int_{x-\zeta}^{x} \tilde{\omega}(x-s)\rho(t,s)ds \right), \text{ for } t \geq 0, x \in \Re \tag{2.3}$$

where $\zeta > 0$ is a constant, $g : \Re_+ \to \Re_+$ is a non-decreasing, bounded function, and $\tilde{\omega} : \Re_+ \to \Re_+$ is a non-increasing function. Even more emphatically, in the era of automated vehicles, the functions $g : \Re_+ \to \Re_+$ and $\tilde{\omega} : \Re_+ \to \Re_+$ may be designed so that the traffic flow behavior of system (2.1), (2.3) has specific characteristics, e.g., so that the equilibrium point gives maximum flow of vehicles and is globally asymptotically stable. It is clear that in such a case the design problem for $g : \Re_+ \to \Re_+$ and $\tilde{\omega} : \Re_+ \to \Re_+$ is strongly reminiscent of the feedback design problem for control systems. Therefore, there is an interest to understand traffic flow models of the form (2.1), (2.3). The first thing that we need to understand is the set of properties that all the functions described above must possess in order to have a well-defined system with solutions that have physical meaning (e.g., $\rho(t,x)$, $v(t,x)$ have to be positive).

In this paper, we study traffic flow models on a ring-road; hence we impose the periodicity condition

$$\rho(t,x+1)=\rho(t,x)\text{, for } t\geq 0, x\in\Re . \tag{2.4}$$

Moreover, we allow the density-speed relation to be of the form

$$v(t,x)=(K(\rho[t]))(x)\text{, for } t\geq 0, x\in\Re \tag{2.5}$$

where $K:Per(\Re)\to C^1(\Re)\cap Per(\Re)$ is a mapping of class $K\in C^0(Per(\Re);Per(\Re))$, for which there exists a constant $v_{\max}>0$ such that the inequality $0\leq (K(\rho))(x)\leq v_{\max}$ holds for all $\rho\in Per(\Re)$. We assume the existence of a non-decreasing function $a:\Re_+\to\Re_+$ such that for every $\rho,\bar{\rho}\in Per(\Re)$ the following inequalities hold:

$$\begin{aligned}&\int_0^1 |\rho(x)-\bar{\rho}(x)||(K(\rho))(x)-(K(\bar{\rho}))(x)|dx\leq a\left(\|\rho\|_\infty+\|\bar{\rho}\|_\infty\right)\int_0^1|\rho(x)-\bar{\rho}(x)|^2dx\\ &\int_0^1 |\rho(x)-\bar{\rho}(x)|\left|\frac{\partial}{\partial x}\left((K(\rho))(x)-(K(\bar{\rho}))(x)\right)\right|dx\leq a\left(\|\rho\|_\infty+\|\bar{\rho}\|_\infty\right)\int_0^1|\rho(x)-\bar{\rho}(x)|^2dx\end{aligned} \tag{2.6}$$

Inequalities (2.6) are technical conditions needed for uniqueness of solutions. It will be shown (in the proof of Theorem 2.4 below) that inequalities (2.6) hold when (2.3) holds, i.e., when $(K(\rho))(x)=f\left(\int_x^{x+\eta}\omega(s-x)\rho(s)ds\right)g\left(\int_{x-\zeta}^{x}\tilde{\omega}(x-s)\rho(s)ds\right)$, under mild assumptions for the functions $f,g$.

Given $\rho_0\in Per(\Re)$ we consider the initial-value problem (2.1), (2.4), (2.5) with initial condition

$$\rho[0]=\rho_0 \tag{2.7}$$

For the statement of our first main result we need the following definition.

**Definition 2.1:** *Suppose that $K\in C^0(Per(\Re);Per(\Re))$ is a mapping for which there exists a constant $v_{\max}>0$ such that the inequality $0\leq (K(\rho))(x)\leq v_{\max}$ holds for all $x\in\Re$, $\rho\in Per(\Re)$. Moreover suppose that there exists a parameterized family of mappings $K_N:\Re^N\to\Re_+$ with parameter the integer $N>2$, constants $L,c,C,S\geq 0$ and non-decreasing functions $\gamma,\Gamma,W:\Re_+\to\Re_+$ with the following property:*

**(P)** *For each $N>2$ there exists a mapping $K_N:\Re^N\to\Re_+$ with $K_N(\rho)\leq v_{\max}$ for all $\rho\in\Re_+^N$, such that the following inequalities hold for all $N>2$, $0<\rho_{\min}\leq\rho_{\max}$, $i=0,...,N-1$ and $\rho,\tilde{\rho}\in\Re_+^N$ with $\rho_{\min}1_N\leq\rho=(\rho_0,...,\rho_{N-1})^\top\leq\rho_{\max}1_N$, $\rho_{\min}1_N\leq\tilde{\rho}\leq\rho_{\max}1_N$:*

$$|K_N(\rho)-K_N(\tilde{\rho})|\leq L|\rho-\tilde{\rho}|_\infty \tag{2.8}$$

$$\left|K_N(\rho)-K_N(\rho^{(1)})-K_N(\tilde{\rho})+K_N(\tilde{\rho}^{(1)})\right|\leq\left(h|\rho-\tilde{\rho}|_\infty+h^2\right)\Gamma(\rho_{\max}) \tag{2.9}$$

$$-ch\left(\rho_{\max}-\rho_0\right)\leq K_N(\rho^{(1)})-K_N(\rho)\leq ch\left(\rho_0-\rho_{\min}\right) \tag{2.10}$$

$$\left|2K_N(\rho^{(1)})-K_N(\rho)-K_N(\rho^{(2)})\right|\leq h^2\left(\gamma(\rho_{\max})+C|y|_\infty\right) \tag{2.11}$$

$$\left|3K_N(\rho^{(1)})+K_N(\rho^{(3)})-K_N(\rho)-3K_N(\rho^{(2)})\right|\leq h^3\left(W(\rho_{\max}+|y|_\infty)+C|\varphi|_\infty\right) \tag{2.12}$$

$$\left|(K(P_N\rho))(ih) - K_N\left(\rho^{(i)}\right)\right| \le hS\rho_{\max} \tag{2.13}$$

*where* $h = 1/N$, $\rho^{(i)} = (\rho_i,...,\rho_{N-1},\rho_0,...,\rho_{i-1})$ *for* $i = 1,...,N-1$, $y = (y_0,...,y_{N-1})^{\top} = h^{-1}(\rho^{(1)} - \rho)$, $\varphi = (\varphi_0,...,\varphi_{N-1})^{\top} = h^{-2}(\rho^{(2)} - 2\rho^{(1)} + \rho)$ *and* $P_N\rho \in Per(\Re)$ *is the function defined by the equations:*

$$(P_N\rho)(x) = h^{-1}\rho_i((i+1)h - x) + \rho_{i+1}h^{-1}(x - ih), \text{ for all } x \in [ih,(i+1)h)\ i = 0,...,N-2 \tag{2.14}$$

$$(P_N\rho)(x) = h^{-1}\rho_{N-1}(1-x) + \rho_0 h^{-1}(x+h-1), (N-1)h \le x \le 1, \text{ for all } x \in [(N-1)h,1] \tag{2.15}$$

*Then we say that the family* $K_N$ *smoothly approximates the mapping* $K \in C^0(Per(\Re);Per(\Re))$.

The class of mappings $K \in C^0(Per(\Re);Per(\Re))$ which can be approximated smoothly by a parameterized family $K_N$ is closed under multiplication and addition. More precisely, we have the following lemma, which is proved in the Appendix.

**Lemma 2.2:** *Suppose that the mappings* $K,G \in C^0(Per(\Re);Per(\Re))$ *can be approximated smoothly by the parameterized families* $K_N, G_N$. *Let* $\lambda \ge 0$ *be any given real number. Then the mappings* $(KG) \in C^0(Per(\Re);Per(\Re))$, $(K+G) \in C^0(Per(\Re);Per(\Re))$, $(\lambda K) \in C^0(Per(\Re);Per(\Re))$ *defined by* $((KG)(\rho))(x) = (K(\rho))(x)(G(\rho))(x)$, $((K+G)(\rho))(x) = (K(\rho))(x) + (G(\rho))(x)$, $((\lambda K)(\rho))(x) = \lambda(K(\rho))(x)$ *for all* $x \in \Re$, $\rho \in Per(\Re)$, *can be approximated smoothly by the parameterized families* $K_N G_N$, $K_N + G_N$, $\lambda K_N$, *respectively.*

Our first main result is an existence and uniqueness result for the initial-value problem (2.1), (2.4), (2.5), (2.7).

**Theorem 2.3:** *Suppose that* $K : Per(\Re) \to C^1(\Re) \cap Per(\Re)$ *is a mapping of class* $K \in C^0(Per(\Re);Per(\Re))$ *for which there exists a constant* $v_{\max} > 0$ *such that the inequality* $0 \le (K(\rho))(x) \le v_{\max}$ *holds for all* $\rho \in Per(\Re)$. *Moreover, suppose that there exists a non-decreasing function* $a : \Re_+ \to \Re_+$ *such that (2.6) holds for every* $\rho, \bar{\rho} \in Per(\Re)$. *Finally, suppose that there exists a parameterized family* $K_N$ *that smoothly approximates the mapping* $K \in C^0(Per(\Re);Per(\Re))$. *Then for every* $\rho_0 \in W^{2,\infty}(\Re) \cap Per(\Re)$ *the initial-value problem (2.1), (2.4), (2.5), (2.7) has a unique solution* $\rho \in C^1(\Re_+ \times \Re)$ *with* $\rho[t] \in W^{2,\infty}(\Re) \cap Per(\Re)$ *for all* $t \ge 0$. *Moreover, the following inequality holds for all* $t \ge 0, x \in \Re$:

$$\min_{x \in [0,1]}\left(\rho_0(x)\right) \le \rho(t,x) \le \max_{x \in [0,1]}\left(\rho_0(x)\right) \tag{2.16}$$

**Remarks: (i)** Theorem 2.3 guarantees the existence of a classical solution for the initial-value problem (2.1), (2.4), (2.5), (2.7). This feature differentiates Theorem 2.3 from other results in the literature (e.g. the results in [3]). Theorem 2.3 shows that the state space for system (2.1), (2.4), (2.5) is the space $W^{2,\infty}(\Re) \cap Per(\Re)$, i.e., if $\rho_0 \in W^{2,\infty}(\Re) \cap Per(\Re)$ then $\rho[t] \in W^{2,\infty}(\Re) \cap Per(\Re)$ for all $t \ge 0$.

**(ii)** The proof of Theorem 2.3 is based on the method of finite differences used in the book [10]. This feature differentiates Theorem 2.3 from other results in the literature where fixed-point theorems are employed.

Theorem 2.3 can be used as a tool for the proof of existence and uniqueness of solutions for system (2.1), (2.3), (2.4). This is achieved by the following theorem, which is the second main result of the paper.

**Theorem 2.4:** *Suppose that* $\eta, \zeta \in (0,1]$, $f, g \in C^3(\Re_+)$, $f'(\rho) \le 0$, $g'(\rho) \ge 0$, $f(\rho) \ge 0$, $g(\rho) \ge 1$ *for all* $\rho \ge 0$. *Moreover, suppose that there exists a constant* $M \ge 0$ *such that*

$$\sup_{\rho \ge 0}\left( \sum_{k=0}^{3} \left| f^{(k)}(\rho) \right| \right) + \sup_{\rho \ge 0}\left( \sum_{k=0}^{3} \left| g^{(k)}(\rho) \right| \right) \le M \tag{2.17}$$

*Finally, suppose that the restrictions of* $\omega, \tilde{\omega} : \Re_+ \to \Re_+$ *on* $[0,\eta]$, $[0,\zeta]$, *respectively, are* $C^1$ *functions with* $\omega'(x) \le 0$ *for* $x \in [0,\eta]$, $\tilde{\omega}'(x) \le 0$ *for* $x \in [0,\zeta]$ *and that* $\omega(x) = 0$ *for* $x > \eta$, $\tilde{\omega}(x) = 0$ *for* $x > \zeta$. *Then for every* $\rho_0 \in W^{2,\infty}(\Re) \cap Per(\Re)$ *the initial-value problem (2.1), (2.4), (2.3), (2.7) has a unique solution* $\rho \in C^1(\Re_+ \times \Re)$ *with* $\rho[t] \in W^{2,\infty}(\Re) \cap Per(\Re)$ *for all* $t \ge 0$. *Moreover, inequality (2.16) holds for all* $t \ge 0, x \in \Re$.

Theorem 2.4 is proved in Section 4 by applying Theorem 2.3 to the case (2.3).

## 3. Controlling Non-Local Traffic Flow Models

The uniform equilibrium points $\rho(x) \equiv \rho > 0$ of model (2.1), (2.2), (2.4) satisfy exactly the same density-flow relation $q = \rho f(\rho)$ of the classical LWR model (the so-called fundamental diagram). This is not true for the uniform equilibrium points $\rho(x) \equiv \rho > 0$ of model (2.1), (2.3), (2.4). For this model, the density-flow relation is given by

$$q = \rho f(\rho) g(\sigma\rho) \tag{3.1}$$

where $\sigma := \int_0^{\zeta} \tilde{\omega}(s) ds$. Since $g(\rho) \ge 1$, relation (3.1) shows that nudging can increase the flow. Moreover, the critical density, i.e., the density for which the flow becomes maximum, changes. This is demonstrated in Fig.1 for the case $f(\rho) = \exp(-\rho)$, $g(s) = (1+k)\frac{\exp(\gamma s)}{k + \exp(\gamma s)}$, $\gamma\sigma = 1$, $k = 1/2$. It may be seen that the flow values are increased, and the critical density is increased as well. It should be noted that Fig.1 is typical for many combinations of functions $f, g$ with the characteristics required by the physics of traffic flow, i.e., $f'(\rho) \le 0$, $g'(\rho) \ge 0$, $f(\rho) > 0$, $g(\rho) \ge 1$, for all $\rho \ge 0$, $\lim_{\rho \to +\infty} \left( f(\rho) \right) = 0$, $\lim_{\rho \to +\infty} \left( g(\rho) \right) < +\infty$.

Theorem 2.4 guarantees that the uniform equilibrium points $\rho(x) \equiv \rho^* > 0$ are neutrally stable in the sup norm of the state. However, Theorem 2.4 says nothing about (local or global) asymptotic stability and convergence to a uniform equilibrium point. Indeed, there are cases where the uniform equilibrium point $\rho(x) \equiv \rho^* > 0$ for model (2.1), (2.2), (2.4) is not locally asymptotically stable, no matter what $f$ is. The following proposition illustrates this point. Its proof is very simple (direct substitution in the equations) and is omitted.

**Proposition 3.1 (Lack of Local Asymptotic Stability for the Model Without Nudging):** *Consider model (2.1), (2.2), (2.4) with* $\omega(x) = \eta^{-1}$ *for* $x \in [0,\eta]$ *and* $\omega(x) = 0$ *for* $x > \eta$, *where* $\eta \in (0,1]$ *is a rational number with* $\eta = \frac{p}{q}$, $p, q > 0$ *integers and* $f \in C^3(\Re_+)$ *is any function with* $f'(\rho) \le 0$, $f(\rho) \ge 0$ *for all* $\rho \ge 0$. *Moreover, suppose that there exists a constant* $M \ge 0$ *such that inequality (2.17) holds with* $g(s) \equiv 1$. *Then, for every* $\rho^* > 0$, $b \in \Re$ *with* $|b| < \rho^*$ *and for every* *integer* $k>0$ *the functions* $\rho(t,x)=\rho^{*}+b\sin\left(2kq\pi\left(x-f(\rho^{*})t\right)\right)$ *for* $t\geq 0, x\in\Re$ *are solutions of (2.1), (2.4), (2.2).*

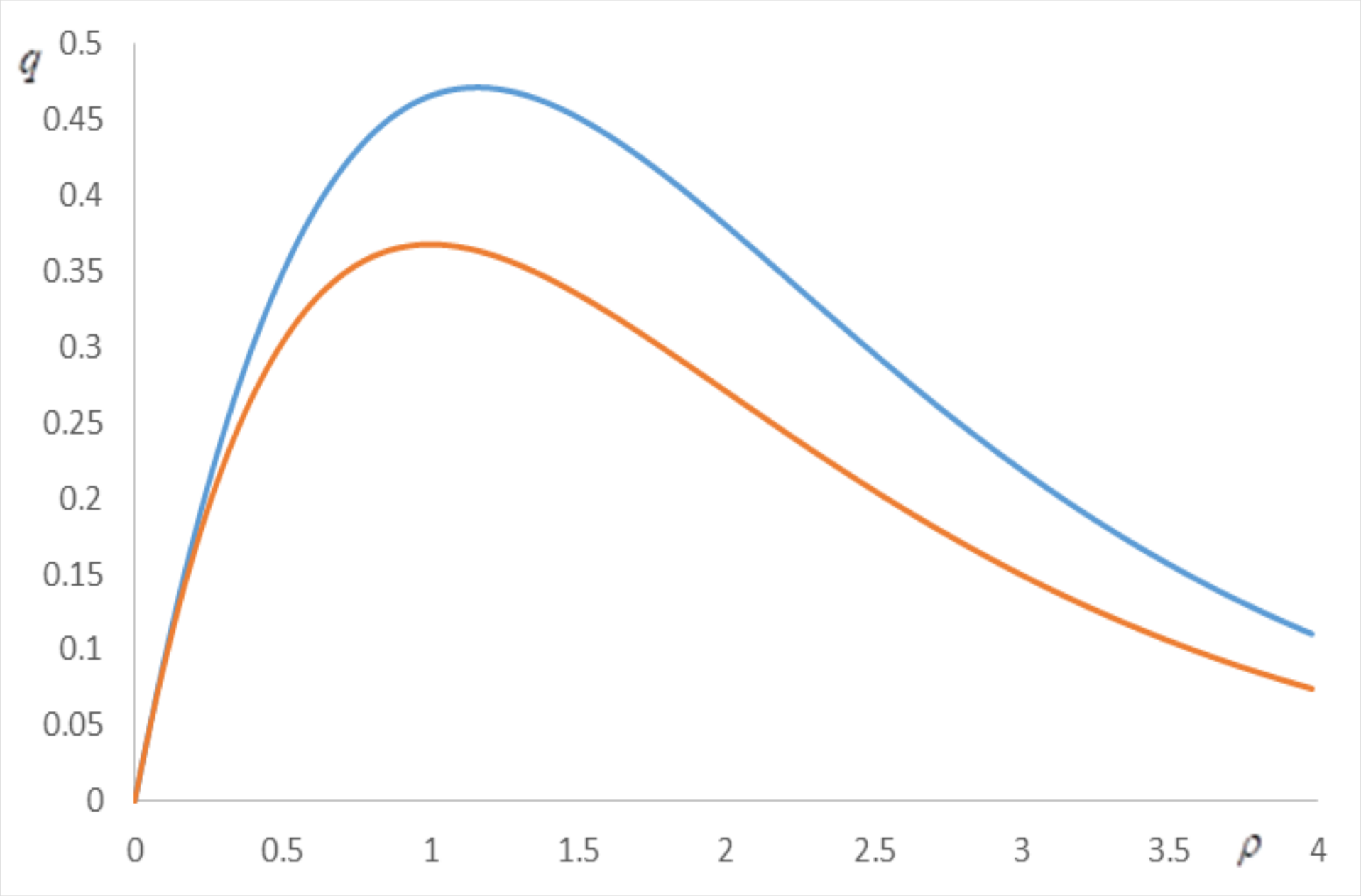

**Fig. 1:** The effect of nudging on the fundamental diagram.
The blue line is the fundamental diagram with nudging and
the red line is the fundamental diagram without nudging ( $g(s)\equiv 1$).

The reader should notice that, if convergence to the uniform equilibrium point $\rho(x)\equiv\rho^{*}>0$ is to be studied, then we should restrict our attention to initial conditions $\rho_0\in W^{2,\infty}(\Re)\cap Per(\Re)$ with $\int_0^1 \rho_0(x)dx=\rho^{*}$, since only for this set of functions we can obtain solutions which converge to the uniform equilibrium point $\rho(x)\equiv\rho^{*}>0$ (notice that $\int_0^1 \rho(t,x)dx=\int_0^1 \rho_0(x)dx$ for all $t\geq 0$ for every solution of (2.1), (2.4), (2.2), (2.7) or any solution of (2.1), (2.4), (2.3), (2.7)).

Proposition 3.1 shows that in the important case $\omega(x)=\eta^{-1}$ for $x\in[0,\eta]$ and $\omega(x)=0$ for $x>\eta$ (case studied in [17]), where $\eta\in(0,1]$ is a rational number, there are initial conditions $\rho_0\in W^{2,\infty}(\Re)\cap Per(\Re)$ with $\int_0^1 \rho_0(x)dx=\rho^{*}$, which are arbitrarily close to the uniform equilibrium point $\rho(x)\equiv\rho^{*}>0$ (in any space $L^{s}(0,1)$ with $s\geq 1$), for which the solution of (2.1), (2.4), (2.2), (2.7) does not converge (in any space $L^{s}(0,1)$ with $s\geq 1$) to the uniform equilibrium point $\rho(x)\equiv\rho^{*}>0$. The solutions that fail to converge to the uniform equilibrium point $\rho(x)\equiv\rho^{*}>0$ are high-frequency density waves which move with constant speed $v(t,x)\equiv f(\rho^{*})$. Since the set of rational numbers is dense within the reals, it follows that we can never be sure about the existence or not of such solutions (no matter how accurate is the measurement of $\eta\in(0,1]$).

In such cases, if properly designed, nudging can improve the stability properties of the system. This is shown by the following result.

**Theorem 3.2 (Local Stabilization by Means of Nudging):** *Consider model (2.1), (2.3), (2.4) with* $\zeta = 1$, $\omega(x) = \eta^{-1}$ *for* $x \in [0,\eta]$ *and* $\omega(x) = 0$ *for* $x > \eta$, $\tilde{\omega}(x) = 1 - x$ *for* $x \in [0,1]$, *where* $\eta \in (0,1]$ *is a constant and* $f, g \in C^3(\Re_+)$ *are any functions with* $f'(\rho) < 0$, $g'(\rho) > 0$, $f(\rho) > 0$, $g(\rho) \geq 1$ *for all* $\rho \geq 0$. *Moreover, suppose that there exists a constant* $M \geq 0$ *such that inequality (2.17) holds. Let* $\rho_0 \in W^{2,\infty}(\Re) \cap Per(\Re)$ *with*

$$F_{\max} g_{\max} - F_{\min} g_{\min} < 2\eta f_{\min} G_{\min} \tag{3.2}$$

*where*

$$F_{\max} := \max\left\{ |f'(s)| : \rho_{\min} \leq s \leq \min\left(\eta^{-1}\rho^*, \rho_{\max}\right) \right\} \tag{3.3}$$

$$F_{\min} := \min\left\{ |f'(s)| : \rho_{\min} \leq s \leq \min\left(\eta^{-1}\rho^*, \rho_{\max}\right) \right\} \tag{3.4}$$

$$f_{\min} := f\left(\min\left(\eta^{-1}\rho^*, \rho_{\max}\right)\right) \tag{3.5}$$

$$g_{\max} := g\left(\frac{1}{2}\min\left(2\rho^* - \rho_{\min}, \rho_{\max}\right)\right) \tag{3.6}$$

$$g_{\min} := g\left(\frac{1}{2}\max\left(2\rho^* - \rho_{\max}, \rho_{\min}\right)\right) \tag{3.7}$$

$$G_{\min} := \min\left\{ g'(s) : \max\left(2\rho^* - \rho_{\max}, \rho_{\min}\right) \leq 2s \leq \min\left(2\rho^* - \rho_{\min}, \rho_{\max}\right) \right\} \tag{3.8}$$

$\rho^* = \int_0^1 \rho_0(s)ds$, $\rho_{\min} := \min_{x\in[0,1]}\left(\rho_0(x)\right)$ *and* $\rho_{\max} := \max_{x\in[0,1]}\left(\rho_0(x)\right)$. *Then there exists a constant* $\bar{c} > 0$ *such that the unique solution* $\rho \in C^1(\Re_+ \times \Re)$ *of the initial-value problem (2.1), (2.4), (2.3), (2.7) satisfies the estimate:*

$$\int_0^1 \left(\rho(t,x) - \rho^*\right)^2 dx \leq \frac{\rho_{\max}}{\rho_{\min}} \exp(-\bar{c}t) \int_0^1 \left(\rho_0(x) - \rho^*\right)^2 dx \text{, for } t \geq 0 \tag{3.9}$$

**Remarks: (i)** When $\rho_{\min} = \rho_{\max} = \rho^*$ then (3.2) holds automatically (by virtue of the fact that $f'(\rho) < 0$, $g'(\rho) > 0$, $f(\rho) > 0$, $g(\rho) \geq 1$ for all $\rho \geq 0$). Due to continuity of $F_{\max}, F_{\min}, f_{\min}, g_{\max}, g_{\min}, G_{\min}$ with respect to $\rho_{\min}, \rho_{\max}, \rho^*$, for every $\rho^* > 0$, there exist $\rho_{\min} < \rho^* < \rho_{\max}$ such that (3.2) holds. This implies the existence of a neighborhood of the uniform equilibrium point $\rho(x) \equiv \rho^* > 0$ in $Per(\Re)$ for which the $L^2(0,1)$ norm of the deviation of the solution from the equilibrium point converges exponentially to zero.

**(ii)** Condition (3.2) is a condition on the maximum deviation $\left\|\rho_0 - \rho^*\right\|_\infty$ of the initial condition from the desired uniform equilibrium point $\rho(x) \equiv \rho^* > 0$. To see this, notice that condition (3.2) takes the following form when $f''(\rho) \geq 0$ and $g''(\rho) \leq 0$ for all $\rho \geq 0$:

$$f'\left(\min\left(\eta^{-1}\rho^*, \rho_{\max}\right)\right) g\left(\frac{1}{2}\max\left(2\rho^* - \rho_{\max}, \rho_{\min}\right)\right) - f'\left(\rho_{\min}\right) g\left(\frac{1}{2}\min\left(2\rho^* - \rho_{\min}, \rho_{\max}\right)\right)$$
$$< 2\eta f\left(\min\left(\eta^{-1}\rho^*, \rho_{\max}\right)\right) g'\left(\frac{1}{2}\min\left(2\rho^* - \rho_{\min}, \rho_{\max}\right)\right)$$

Therefore, for the case $f(\rho) = A\exp(-b\rho)$, $g(s) = \frac{1+as}{1+\gamma s}$, with $a > \gamma$, $A, b > 0$ and $\eta \leq 1/2$, condition (3.2) becomes:

$$b\exp\left(2b\left\|\rho_0 - \rho^*\right\|_\infty\right) \frac{c_1 + a\left\|\rho_0 - \rho^*\right\|_\infty}{c_2 + \gamma\left\|\rho_0 - \rho^*\right\|_\infty} < b\frac{c_1 - a\left\|\rho_0 - \rho^*\right\|_\infty}{c_2 - \gamma\left\|\rho_0 - \rho^*\right\|_\infty} + \frac{8\eta(a-\gamma)}{\left(c_2 + \gamma\left\|\rho_0 - \rho^*\right\|_\infty\right)^2}$$

where $c_1 = 2 + a\rho^* > c_2 = 2 + \gamma\rho^*$. The above inequality provides a bound $R > 0$ on $\left\|\rho_0 - \rho^*\right\|_\infty$ such that (3.2) (and consequently the exponential stability estimate (3.9)) holds for all $\rho_0 \in W^{2,\infty}(\Re) \cap Per(\Re)$ with $\left\|\rho_0 - \rho^*\right\|_\infty < R$. However, it should be remarked at this point that simulations indicate that condition (3.2) provides a conservative estimate of the bound $R > 0$ (see Section 5).

**(iii)** The proof of Theorem 3.2 shows that $\bar{c} := \eta^{-1}\rho_{\min}\left(2\eta f_{\min}G_{\min} - F_{\max}g_{\max} + F_{\min}g_{\min}\right)$.

**(iv)** The nudging term with $\zeta = 1$, $\tilde{\omega}(x) = 1 - x$ for $x \in [0,1]$, depends on the whole density profile of the ring-road. Such a term has no meaning when the vehicles are driven by human drivers. However, when automated vehicles are present in a highway, then such a term can be implemented by providing continuously information for the density profile to each vehicle. In such a case, the effect of nudging is not only the increase of the flow, but also the elimination of the well-known stop-and-go waves (see [1]).

**(v)** The proof of Theorem 3.2 makes use of estimate (2.16) and the functional $V(\rho) = \int_0^1 \left( \rho(x) \ln\left( \frac{\rho(x)}{\rho^*} \right) + \rho^* - \rho(x) \right) dx$. This functional, defined on the set of functions $\rho \in Per(\Re)$ with $\rho^* = \int_0^1 \rho(x)dx$, is a non-coercive Control Lyapunov Functional for the control system (2.1), (2.4) with $v[t] \in Per(\Re) \cap C^1(\Re)$ as input. Indeed, for classical solutions of (2.1), (2.4) we get $\frac{d}{dt}V(\rho[t]) = \int_0^1 \left(\rho^* - \rho(t,x)\right)\frac{\partial v}{\partial x}(t,x)dx$. For non-coercive Lyapunov functionals, the reader can consult [9, 20]. It is possible that the use of other Lyapunov functionals can give less demanding conditions than (3.2) for exponential convergence to the uniform equilibrium point $\rho(x) \equiv \rho^* > 0$.

## 4. Illustrative Examples

In this section we present some numerical examples that demonstrate the advantages and the stabilizing effects of nudging in comparison with the "look-ahead" model (2.1), (2.2) and the LWR model (2.1) with $v(t,x) = f(\rho(t,x))$. Hence, we consider the three models displayed in Table 1 on a ring road (2.4) with initial density profile

$$\rho_0(x) = \begin{cases} 2.35, & 0.5 \le x \le 0.75 \\ 0.55 & \text{else} \end{cases} \tag{4.1}$$

| | Model 1 | Model 2 | Model 3 |
|---|---|---|---|
| Velocity $v(t,x)$ | LWR with $v(t,x) = f(\rho(t,x))$ | Model (2.1) with $v(t,x)$ given by (2.2) | Model (2.1) with $v(t,x)$ given by (2.3) |
| $f(\cdot)$ | $f(x) = 0.96\exp(-x)$ | $f(x) = 0.96\exp(-x)$ | $f(x) = 0.75\exp(-x)$ |
| $g(\cdot)$ | N/A | N/A | $g(x) = 1.6\dfrac{\exp\left(\dfrac{1.8}{\zeta(2-\zeta)}x\right)}{0.6 + \exp\left(\dfrac{1.8}{\zeta(2-\zeta)}x\right)}$ |
| $\omega(\cdot)$ | N/A | $\eta^{-1}$ | $\eta^{-1}$ |
| $\tilde{\omega}(\cdot)$ | N/A | N/A | $1 - x$ |

**Table 1.** The three models of the simulation examples.

Notice that the initial condition corresponds to a road with congestion belt at $[0.5, 0.75]$ with the same uniform equilibrium for all three models given by $\rho^* = \int_0^1 \rho_0(x)dx = 1$. It should also be noticed that while all models have the same equilibrium, they do not possess the same fundamental diagram and the free road speed for Model 1 and Model 2 is 60% higher than the free road speed in Model 3. We have used two different values for the upstream horizon: $\zeta = 1$ and $\zeta = 0.154$. While Theorem 3.2 guarantees local exponential stabilization for $\zeta = 1$, it is important for implementation purposes to consider small values for the upstream horizon (which do not require knowledge of the whole density profile of the ring road). The value of the downstream horizon in all experiments was set $\eta = 0.1$.

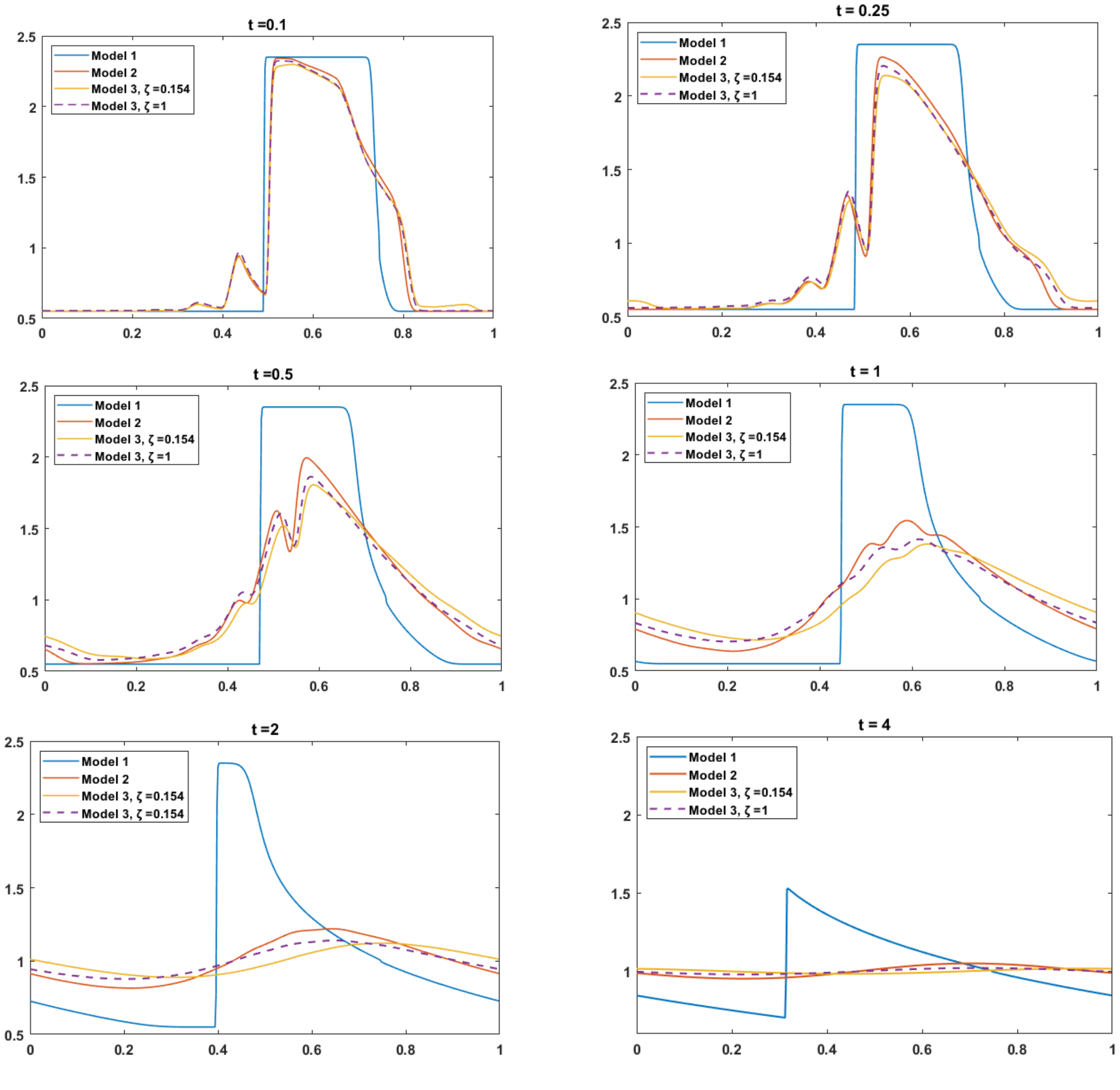

**Fig 2.** Density profiles of Model 1, Model 2, Model 3 with $\zeta = 1$ and Model 3 with $\zeta = 0.154$.

While the numerical results for the LWR model are obtained by means of the Godunov numerical scheme, for the non-local PDEs we have used the numerical scheme (5.3)-(5.10), (5.60) with $h = 1/500$, $\lambda = 0.25$. Figure 2 shows the density profiles at different times. Notice that the rate of convergence of Model 3 for both values $\zeta = 1$ and $\zeta = 0.154$ is faster compared to the other models. This feature can also be verified in Figure 3, which depicts the evolution of the $L^2$ norm of the deviation from the equilibrium $\left\| \rho[t] - \rho^* \right\|_2 = \left( \int_0^1 \left( \rho(t,x) - \rho^* \right)^2 dx \right)^{1/2}$. Figure 3 also shows the evolution of the logarithm of the $L^2$ norm $\left\| \rho[t] - \rho^* \right\|_2$ indicating exponential convergence.

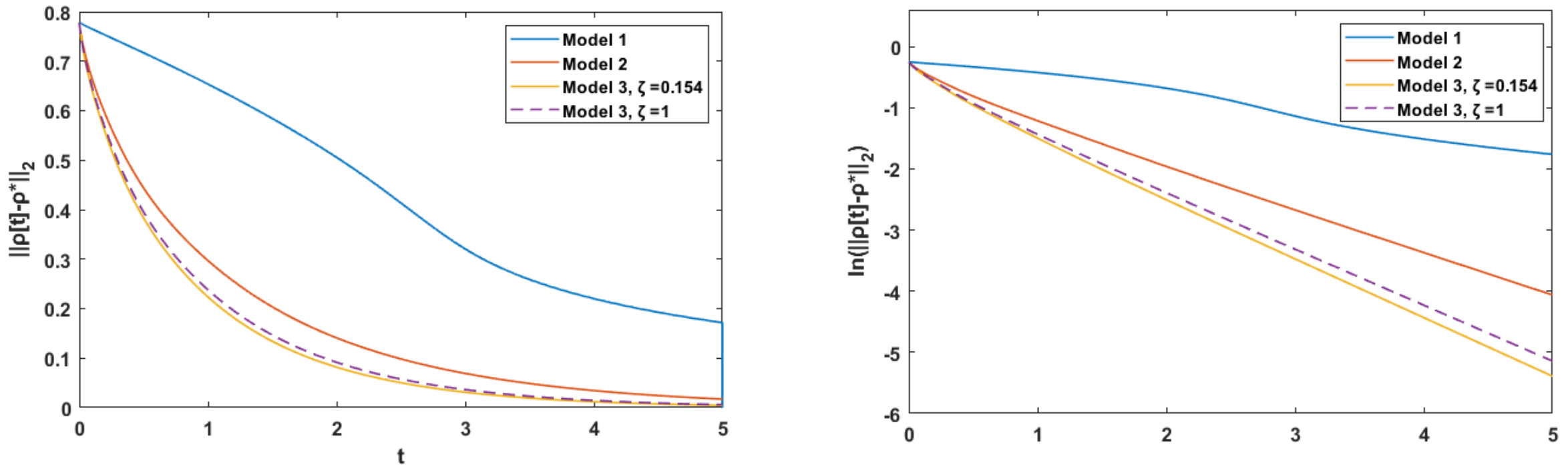

**Fig. 3** Evolution of the $L^2$ norm $\left\| \rho[t] - \rho^* \right\|_2 = \left( \int_0^1 \left( \rho(t,x) - \rho^* \right)^2 dx \right)^{1/2}$ and its logarithm for Model 1, Model 2, Model 3 with $\zeta = 1$ and Model 3 with $\zeta = 0.154$.

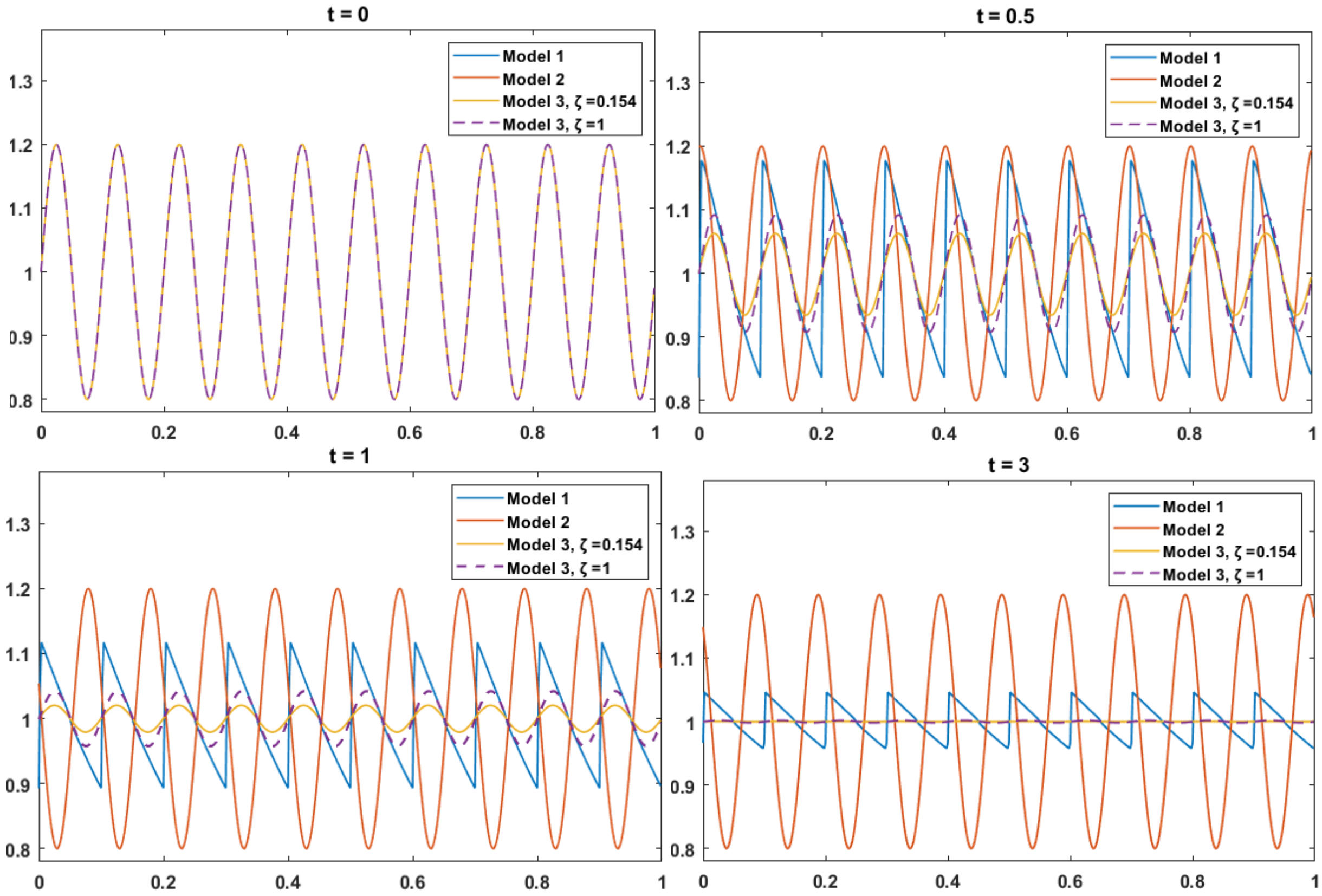

**Fig. 4** Density profiles at various time instants for the initial condition $\rho_0(x) = 1 + 0.2\sin(20\pi x)$.

Finally, we consider again Model 1, Model 2 and Model 3 shown in Table 1, with initial condition $\rho_0(x) = 1 + 0.2\sin(20\pi x)$. In this case Proposition 3.1 holds and the solution for Model 2 is given by the formula $\rho(t,x) = \rho^* + b\sin\left(2kq\pi\left(x - f(\rho^*)t\right)\right)$ with $\rho^* = 1$, $b = 0.2$, $k = 1$ and $q = 10$. The density waves for Model 2 and the density profiles for all other models are displayed in Figure 4. Figure 5 depicts again the evolution of the $L^2$ norm of the deviation from the equilibrium $\left\|\rho[t]-\rho^*\right\|_2=\left(\int_0^1\left(\rho(t,x)-\rho^*\right)^2dx\right)^{1/2}$. Figure 5 also shows the evolution of the logarithm of the $L^2$ norm $\left\|\rho[t]-\rho^*\right\|_2$ indicating exponential convergence for Model 3 with $\zeta=1$ and Model 3 with $\zeta=0.154$.

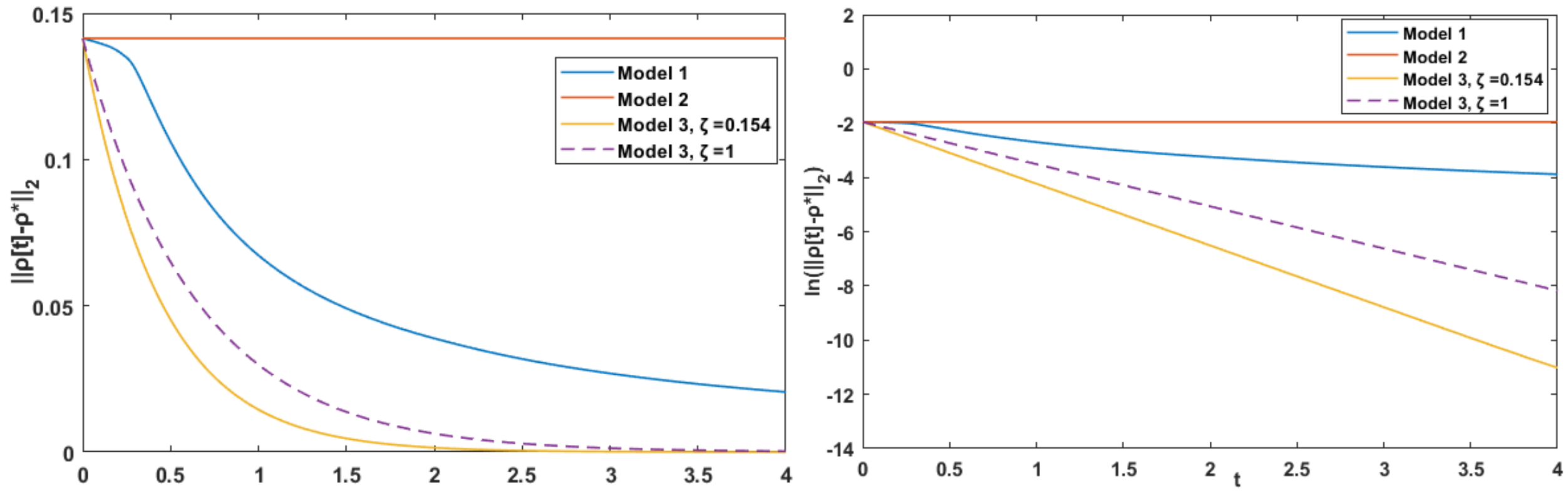

**Fig. 5** Evolution of the $L^2$ norm $\left\|\rho[t]-\rho^*\right\|_2=\left(\int_0^1\left(\rho(t,x)-\rho^*\right)^2dx\right)^{1/2}$ and its logarithm for Model 1, Model 2, Model 3 with $\zeta=1$ and Model 3 with $\zeta=0.154$. Initial condition $\rho_0(x)=1+0.2\sin(20\pi x)$.

## 5. Proofs of Main Results

The proof of Theorem 2.3 requires a technical result, whose proof is simple and is omitted.

**Lemma 5.1:** *Suppose that there exist constants $a>0, b\geq 0$ such that the sequence $\{x(k)\geq 0\}_{k=0}^{\infty}$ satisfies the inequality:*

$$x(k+1)\leq(1+a)x(k)+b\text{, for all } k=0,1,...,m-1 \tag{5.1}$$

*Then the following estimate holds:*

$$x(k)\leq\exp(ka)\left(x(0)+\frac{b}{a}\right)\text{, for all } k=0,1,...,m \tag{5.2}$$

We are now ready to give the proof of Theorem 2.3.

**Proof of Theorem 2.3:** Let $N>2$ and $T>0$ be given. Consider the parameterized infinite-dimensional, discrete-time system

$$\rho_i((k+1)\delta)=\left(1-\lambda v_{i+1}(k\delta)\right)\rho_i(k\delta)+\lambda v_i(k\delta)\rho_{i-1}(k\delta)\text{, for } i=0,\pm1,\pm2,...,\ k=0,1,...,m \tag{5.3}$$

$$\rho_i(0):=\rho_0(ih)\text{, for } i=0,\pm1,\pm2,... \tag{5.4}$$

where

$$h:=1/N,\ \delta:=\lambda h \tag{5.5}$$

$$v_i(k\delta)=K_N((Q\rho(k\delta))^{(i)})\text{, for } i=0,\pm1,\pm2,...,\ k=0,1,...,m+1 \tag{5.6}$$

$$\lambda:=\frac{T}{[T(v_{\max}+c\rho_{\max})]+1} \tag{5.7}$$

$$m:=N\left([T(v_{\max}+c\rho_{\max})]+1\right) \tag{5.8}$$

$$Q\rho(k\delta) = (\rho_0(k\delta),...,\rho_{N-1}(k\delta))^\top \in \Re^N \ , \ k = 0,1,...,m \tag{5.9}$$

$$\rho_{\min} := \min_{x\in[0,1]}\left(\rho_0(x)\right), \ \rho_{\max} := \max_{x\in[0,1]}\left(\rho_0(x)\right) \tag{5.10}$$

Notice that the above definitions guarantee that

$$T = m\delta \ , \tag{5.11}$$

$$\lambda\, v_{\max} + \delta\, c\rho_{\max} \leq 1 . \tag{5.12}$$

We next prove by induction that

$$\rho_{\min} \leq \rho_i(k\delta) \leq \rho_{\max} \text{, for } i = 0,\pm1,\pm2,... , \ k = 0,1,...,m+1 \tag{5.13}$$

Indeed, by virtue of definitions (5.4), (5.10), it follows that (5.13) holds for $k = 0$. Suppose that (5.13) holds for some $k = 0,1,...,m$. Definition (5.9) and property (2.10) imply that

$$hc\left(\rho_i(k\delta) - \rho_{\min}\right) \geq v_{i+1}(k\delta) - v_i(k\delta) \geq -hc\left(\rho_{\max} - \rho_i(k\delta)\right) \text{, for } i = 0,\pm1,\pm2,... \tag{5.14}$$

Using (5.3), (5.5), (5.12), (5.13), (5.14) and the fact that $0 \leq v_i(k\delta) \leq v_{\max}$ for all $i = 0,\pm1,\pm2,...$, $k = 0,1,...,m+1$ (a consequence of (5.6)), we get:

$$\begin{aligned}
&\rho_i((k+1)\delta) = (1-\lambda v_i(k\delta))\rho_i(k\delta) + \lambda v_i(k\delta)\rho_{i-1}(k\delta) - \lambda(v_{i+1}(k\delta) - v_i(k\delta))\rho_i(k\delta)\\
&\leq (1-\lambda v_i(k\delta))\rho_i(k\delta) + \lambda v_i(k\delta)\rho_{i-1}(k\delta) + c\delta\,\rho_i(k\delta)\left(\rho_{\max} - \rho_i(k\delta)\right)\\
&\leq (1-\lambda v_i(k\delta) - c\delta\,\rho_i(k\delta))\rho_i(k\delta) + \left(\lambda v_i(k\delta) + c\delta\,\rho_i(k\delta)\right)\rho_{\max} \leq \rho_{\max}
\end{aligned}$$

$$\begin{aligned}
&\rho_i((k+1)\delta) = (1-\lambda v_i(k\delta))\rho_i(k\delta) + \lambda v_i(k\delta)\rho_{i-1}(k\delta) - \lambda(v_{i+1}(k\delta) - v_i(k\delta))\rho_i(k\delta)\\
&\geq (1-\lambda v_i(k\delta))\rho_i(k\delta) + \lambda v_i(k\delta)\rho_{i-1}(k\delta) - c\delta\,\rho_i(k\delta)\left(\rho_i(k\delta) - \rho_{\min}\right)\\
&\geq (1-\lambda v_i(k\delta) - c\delta\,\rho_i(k\delta))\rho_i(k\delta) + \left(\lambda v_i(k\delta) + c\delta\,\rho_i(k\delta)\right)\rho_{\min} \geq \rho_{\min}
\end{aligned}$$

Consequently, (5.13) holds for $k+1$ and this completes the induction step.

Next define

$$y_i(k\delta) := h^{-1}\left(\rho_{i+1}(k\delta) - \rho_i(k\delta)\right) \text{, for } i = 0,\pm1,\pm2,... , \ k = 0,1,...,m \tag{5.15}$$

$$g_i(k\delta) = \delta^{-1}\left(\rho_i((k+1)\delta) - \rho_i(k\delta)\right) \text{, for } i = 0,\pm1,\pm2,... , \ k = 0,1,...,m \tag{5.16}$$

and notice that (5.4) and definition (5.15) imply that $\max_{j=...-1,0,1,...}\left(\left|y_j(0)\right|\right) \leq \left\|\rho_0'\right\|_\infty$. Using (5.3) and definitions (5.15), (5.16), we obtain the following equations:

$$y_i((k+1)\delta) = \left(1-\lambda v_{i+2}(k\delta)\right)y_i(k\delta) + \lambda v_i(k\delta)y_{i-1}(k\delta) + \lambda h^{-1}\left(2v_{i+1}(k\delta) - v_{i+2}(k\delta) - v_i(k\delta)\right)\rho_i(k\delta)$$
$$\text{for } i = 0,\pm1,\pm2,... , \ k = 0,1,...,m-1 \tag{5.17}$$

$$g_i(k\delta) = h^{-1}\left(v_i(k\delta) - v_{i+1}(k\delta)\right)\rho_i(k\delta) - v_i(k\delta)y_{i-1}(k\delta)$$
$$\text{for } i = 0,\pm1,\pm2,... , \ k = 0,1,...,m \tag{5.18}$$

Using (2.11), (5.6) and (5.13), we obtain the following estimate for $i = 0,\pm1,\pm2,...$, $k = 0,1,...,m$:

$$\left|2v_{i+1}(k\delta)-v_i(k\delta)-v_{i+2}(k\delta)\right|\le h^2\gamma\left(\rho_{\max}\right)+h^2C\max_{j=\dots,-1,0,1,\dots}\left(\left|y_j(k\delta)\right|\right) \tag{5.19}$$

It follows from (5.5), (5.17), (5.19), (5.12), (5.13) and (2.10) (which in conjunction with (5.6), (5.13) implies that $\left|v_{i+1}(k\delta)-v_i(k\delta)\right|\le hc\rho_{\max}$ for all $i=0,\pm1,\pm2,\dots$, $k=0,1,\dots,m+1$) that the following inequality holds for $k=0,1,\dots,m-1$:

$$\max_{j=\dots,-1,0,1,\dots}\left(\left|y_j((k+1)\delta)\right|\right)\le\left(1+\delta(2c+C)\rho_{\max}\right)\max_{j=\dots,-1,0,1,\dots}\left(\left|y_j(k\delta)\right|\right)+\delta\rho_{\max}\gamma(\rho_{\max}) \tag{5.20}$$

Using Lemma 2, it follows that inequality (5.20) in conjunction with (5.11) and the fact that $\max_{j=\dots-1,0,1,\dots}\left(\left|y_j(0)\right|\right)\le\left\|\rho_0'\right\|_\infty$, implies the following estimate for $k=0,1,\dots,m$:

$$\max_{j=\dots,-1,0,1,\dots}\left(\left|y_j(k\delta)\right|\right)\le Y:=\exp\left(T(2c+C)\rho_{\max}\right)\left(\left\|\rho_0'\right\|_\infty+(2c+C)^{-1}\gamma(\rho_{\max})\right) \tag{5.21}$$

Moreover, we obtain from (2.8), (5.13), (5.16), (5.18), (5.21) and the facts that $\left|v_{i+1}(k\delta)-v_i(k\delta)\right|\le hc\rho_{\max}$, $v_i(k\delta)\le v_{\max}$ for all $i=0,\pm1,\pm2,\dots$:

$$\max_{i=\dots,-1,0,1,\dots}\left(\left|g_i(k\delta)\right|\right)\le G:=c\rho_{\max}^2+v_{\max}Y\text{, for } k=0,1,\dots,m \tag{5.22}$$

$$\left|v_i((k+1)\delta)-v_i(k\delta)\right|\le GL\delta\text{, for } i=0,\pm1,\pm2,\dots,\ k=0,1,\dots,m \tag{5.23}$$

We define the functions $\rho(t,x;N)$, $v(t,x;N)$ for $(t,x)\in[0,T]\times\Re$ and for every integer $N>2$ (recall that $h=N^{-1}$, $\delta=\lambda h$, $m\delta=T$):

$$\begin{aligned}\rho(k\delta,x;N)&=\left(i+1-xN\right)\rho_i(k\delta)+\left(xN-i\right)\rho_{i+1}(k\delta),\\ v(k\delta,x;N)&=\left(i+1-xN\right)v_i(k\delta)+\left(xN-i\right)v_{i+1}(k\delta)\\ &\text{with } i=[xN]\text{, for } x\in\Re,\ k=0,\dots,m,\end{aligned} \tag{5.24}$$

$$\begin{aligned}\rho(t,x;N)&=\left(k+1-\lambda^{-1}tN\right)\rho(k\delta,x;N)+\left(\lambda^{-1}tN-k\right)\rho((k+1)\delta,x;N)\\ v(t,x;N)&=\left(k+1-\lambda^{-1}tN\right)v(k\delta,x;N)+\left(\lambda^{-1}tN-k\right)v((k+1)\delta,x;N)\\ &\text{with } k=[\lambda^{-1}tN]\text{ for } x\in\Re,\ t\in[0,T).\end{aligned} \tag{5.25}$$

It follows from (5.3), (5.4) and (5.6) that $\rho_{i+N}(k\delta)=\rho_i(k\delta)$, $v_{i+N}(k\delta)=v_i(k\delta)$ for all $i=0,\pm1,\pm2,\dots$, $k=0,1,\dots,m$. Therefore, definitions (5.24), (5.25) imply that for each $t\in[0,T]$ the functions $\rho(t,\cdot;N)$, $v(t,\cdot;N)$ are periodic with period 1. Estimate (5.13), definition (5.6) and the fact that $0\le K_N(\rho)\le v_{\max}$ for all $\rho\in\Re_+^N$, in conjunction with definitions (5.24), (5.25) imply that the following estimates hold for every integer $N>2$:

$$\rho_{\min}\le\rho(t,x;N)\le\rho_{\max}\text{, for } (t,x)\in[0,T]\times\Re. \tag{5.26}$$

$$0\le v(t,x;N)\le v_{\max}\text{, for } (t,x)\in[0,T]\times\Re. \tag{5.27}$$

Definitions (5.15), (5.16) in conjunction with estimates (5.21), (5.22), (5.23) and the fact that $\left|v_{i+1}(k\delta)-v_i(k\delta)\right|\le hc\rho_{\max}$ for all $i=0,\pm1,\pm2,\dots$, imply the existence of a constant $\bar{L}>0$ independent of $N>2$ for which the following estimate holds for all $i,j=0,\pm1,\pm2,\dots$, $k,l=0,1,\dots,m$:

$$\left|\rho_i(k\delta)-\rho_j(l\delta)\right|+\left|v_i(k\delta)-v_j(l\delta)\right|\le\bar{L}(h\left|i-j\right|+\delta\left|k-l\right|) \tag{5.28}$$

Estimate (5.28) in conjunction with definitions (5.24), (5.25) implies that there exists a constant $L_1 > 0$ (independent of $N > 2$) such that the following Lipschitz inequality holds for every integer $N > 2$:

$$|\rho(t,x;N) - \rho(\tau,z;N)| + |v(t,x;N) - v(\tau,z;N)| \le L_1\left(|x-z| + |t-\tau|\right),$$
$$\text{for all } t,\tau \in [0,T] \text{ and } x,z \in \Re. \tag{5.29}$$

It follows from (5.26), (5.27), (5.29) that the sequences of functions $\{v(\cdot;N)\}_{N=3}^{\infty}$, $\{\rho(\cdot;N)\}_{N=3}^{\infty}$ are uniformly bounded and equicontinuous. Therefore, compactness of $[0,T]\times[0,1]$ and the Arzela-Ascoli theorem implies that there exist Lipschitz functions $v:[0,T]\times[0,1]\to\Re$, $\rho:[0,T]\times[0,1]\to\Re$ and subsequences $\{v(\cdot;N_q)\}_{q=1}^{\infty}$, $\{\rho(\cdot;N_q)\}_{q=1}^{\infty}$ for an increasing index sequence $\{N_q\}_{q=1}^{\infty}$, which converge uniformly on $[0,T]\times[0,1]$ to $v,\rho$, respectively. Moreover, the functions $v,\rho$ satisfy the same bounds with $v(\cdot;N),\rho(\cdot;N)$, i.e., $\rho_{\min} \le \rho(t,x) \le \rho_{\max}$, $0 \le v(t,x) \le v_{\max}$, for $(t,x)\in[0,T]\times[0,1]$.

Since the functions $\rho(t,\cdot;N)$, $v(t,\cdot;N)$ are periodic with period 1, it follows that $\rho(t,1)=\rho(t,0)$ and $v(t,1)=v(t,0)$ for all $t\in[0,T]$. Therefore, the subsequences $\{v(\cdot;N_q)\}_{q=1}^{\infty}$, $\{\rho(\cdot;N_q)\}_{q=1}^{\infty}$ converge uniformly on $[0,T]\times\Re$ to the periodic extensions with respect to $x$ (with period 1) of $v,\rho$, respectively. We will denote by $v,\rho$ the periodic extensions with respect to $x$ (with period 1) of $v,\rho$ (a slight abuse of notation).

We show next that (2.5) holds for $(t,x)\in[0,T]\times\Re$. By virtue of (2.13), (5.13) and (5.6) the following inequality holds for $i=0,\pm1,\pm2,...$, $k=0,1,...,m$:

$$|(K(P_N Q\rho(k\delta)))(ih) - v_i(k\delta)| \le hS\rho_{\max} \tag{5.30}$$

Definitions (2.14), (2.15), (5.9), (5.24), (5.25) imply that

$$(P_N Q\rho(k\delta))(x) = \rho(k\delta,x;N), \text{ for all } k=0,1,...,m \text{ and } x\in\Re,\ N>2 \tag{5.31}$$

It follows from (5.30), (5.24), (5.25) that the following equality holds for $i=0,\pm1,\pm2,...$, $k=0,1,...,m$:

$$|(K(\rho(k\delta,\cdot;N)))(ih) - v(k\delta,ih;N)| \le hS\rho_{\max} \tag{5.32}$$

Let $(t,x)\in[0,T]\times\Re$ be given (arbitrary). Let $k=[\lambda^{-1}tN]$, $i=[xN]$ and notice that $|x-ih|\le h=N^{-1}$, $|t-k\delta|\le\delta=\lambda N^{-1}$. We get:

$$\begin{aligned}
|(K(\rho[t]))(x) - v(t,x)| \le\ & |(K(\rho[t]))(x) - (K(\rho[t]))(ih)| \\
&+ |(K(\rho[t]))(ih) - (K(\rho[k\delta]))(ih)| \\
&+ |(K(\rho[k\delta]))(ih) - (K(\rho(k\delta,\cdot;N)))(ih)| \\
&+ |(K(\rho(k\delta,\cdot;N)))(ih) - v(k\delta,ih;N)| \\
&+ |v(k\delta,ih;N) - v(k\delta,ih)| \\
&+ |v(k\delta,ih) - v(t,x)|
\end{aligned}$$

Since

- the subsequences $\{v(\cdot;N_q)\}_{q=1}^{\infty}$, $\{\rho(\cdot;N_q)\}_{q=1}^{\infty}$ converge uniformly on $[0,T]\times\Re$ to $v,\rho$, respectively,
- $K(\rho[t])\in C^1(\Re)\cap Per(\Re)$, which implies the existence of a constant $M>0$ such that $|(K(\rho[t]))(x) - (K(\rho[t]))(y)| \le M|x-y|$,
- $K\in C^0(Per(\Re);Per(\Re))$,
- $v,\rho$ are Lipschitz functions,

it follows from (5.32) that all terms in the right hand side of the above inequality can become arbitrarily small for sufficiently large $q$. Therefore, (2.5) holds for $(t,x)\in[0,T]\times\Re$.

Next define

$$\varphi_i(k\delta)=h^{-1}\left(y_{i+1}(k\delta)-y_i(k\delta)\right),\text{ for } i=0,\pm1,\pm2,...,\ k=0,1,...,m \tag{5.33}$$

$$\psi_i(k\delta)=\delta^{-1}\left(y_i((k+1)\delta)-y_i(k\delta)\right),\text{ for } i=0,\pm1,\pm2,...,\ k=0,1,...,m-1 \tag{5.34}$$

$$\eta_i(k\delta)=h^{-1}\left(g_{i+1}(k\delta)-g_i(k\delta)\right),\text{ for } i=0,\pm1,\pm2,...,\ k=0,1,...,m \tag{5.35}$$

$$\mu_i(k\delta)=\delta^{-1}\left(g_i((k+1)\delta)-g_i(k\delta)\right),\text{ for } i=0,\pm1,\pm2,...,\ k=0,1,...,m-1 \tag{5.36}$$

$$\theta_i(k\delta)=h^{-1}\left(v_{i+1}(k\delta)-v_i(k\delta)\right),\text{ for } i=0,\pm1,\pm2,...,\ k=0,1,...,m \tag{5.37}$$

and notice that (5.4) and definitions (5.15), (5.33) imply that $\max_{j=...-1,0,1,...}\left(\left|\varphi_j(0)\right|\right)\le\left\|\rho_0''\right\|_\infty$. Using definitions (5.33), (5.34), (5.35), (5.36), (5.15), (5.16), we get:

$$\begin{aligned}&\varphi_i((k+1)\delta)=\left(1-\lambda v_{i+1}(k\delta)\right)\varphi_i(k\delta)+\lambda v_{i+1}(k\delta)\varphi_{i-1}(k\delta)\\&+\lambda h^{-2}\left(3v_{i+2}(k\delta)-v_{i+3}(k\delta)-3v_{i+1}(k\delta)+v_i(k\delta)\right)\rho_i(k\delta)\\&+\lambda h^{-1}\left(2v_{i+2}(k\delta)-v_{i+3}(k\delta)-v_{i+1}(k\delta)\right)(y_i(k\delta)+y_{i-1}(k\delta))\\&+\lambda h^{-1}\left(2v_{i+1}(k\delta)-v_i(k\delta)-v_{i+2}(k\delta)\right)y_{i-1}(k\delta)+\lambda\left(v_{i+1}(k\delta)-v_{i+2}(k\delta)\right)\varphi_i(k\delta)\\&+\lambda\left(v_{i+2}(k\delta)-v_{i+3}(k\delta)\right)(\varphi_i(k\delta)+\varphi_{i-1}(k\delta))\\&\qquad\text{for } i=0,\pm1,\pm2,...,\ k=0,1,...,m-1\end{aligned} \tag{5.38}$$

$$\begin{aligned}&\psi_i(k\delta)=-v_{i+2}(k\delta)\varphi_{i-1}(k\delta)+h^{-1}\left(v_i(k\delta)-v_{i+2}(k\delta)\right)y_{i-1}(k\delta)\\&+h^{-2}\left(2v_{i+1}(k\delta)-v_{i+2}(k\delta)-v_i(k\delta)\right)\rho_i(k\delta)\\&\qquad\text{for } i=0,\pm1,\pm2,...,\ k=0,1,...,m-1\end{aligned} \tag{5.39}$$

$$\begin{aligned}&\eta_i(k\delta)=h^{-2}\left(2v_{i+1}(k\delta)-v_{i+2}(k\delta)-v_i(k\delta)\right)\rho_i(k\delta)-v_{i+1}(k\delta)\varphi_{i-1}(k\delta)\\&+h^{-1}\left(v_{i+1}(k\delta)-v_{i+2}(k\delta)\right)y_i(k\delta)+h^{-1}\left(v_i(k\delta)-v_{i+1}(k\delta)\right)y_{i-1}(k\delta)\\&\qquad\text{for } i=0,\pm1,\pm2,...,\ k=0,1,...,m\end{aligned} \tag{5.40}$$

$$\begin{aligned}&\mu_i(k\delta)=h^{-1}\delta^{-1}\left(v_i((k+1)\delta)-v_i(k\delta)+v_{i+1}(k\delta)-v_{i+1}((k+1)\delta)\right)\rho_i(k\delta)\\&+\delta^{-1}\left(v_i((k+1)\delta)-v_{i+1}((k+1)\delta)\right)\lambda g_i(k\delta)\\&-\delta^{-1}\left(v_i((k+1)\delta)-v_i(k\delta)\right)y_{i-1}((k+1)\delta)-v_i(k\delta)\psi_{i-1}(k\delta)\\&\qquad\text{for } i=0,\pm1,\pm2,...,\ k=0,1,...,m-1\end{aligned} \tag{5.41}$$

Using (2.11), (2.12), (5.5), (5.6), (5.12), (5.13), (5.21) and the facts that $0\le v_i(k\delta)\le v_{\max}$, $\left|v_{i+1}(k\delta)-v_i(k\delta)\right|\le hc\rho_{\max}$ for all $i=0,\pm1,\pm2,...$, $k=0,1,...,m$, we get for $i=0,\pm1,\pm2,...$, $k=0,1,...,m-1$:

$$\begin{aligned}&\max_{j=...,-1,0,1,...}\left(\left|\varphi_j((k+1)\delta)\right|\right)\le\left(1+\delta\left(3c+C\right)\rho_{\max}\right)\max_{j=...,-1,0,1,...}\left(\left|\varphi_j(k\delta)\right|\right)\\&+\delta\left(\rho_{\max}W\left(\rho_{\max}+Y\right)+3Y\gamma(\rho_{\max})+3CY^2\right)\end{aligned} \tag{5.42}$$

Using Lemma 2, it follows that inequality (5.42) in conjunction with (5.11) and the fact that $\max_{j=...-1,0,1,...}\left(\left|\varphi_j(0)\right|\right) \le \|\rho_0''\|_\infty$, gives the following estimate for $k = 0,1,...,m$:

$$\max_{j=...,-1,0,1,...}\left(\left|\varphi_j(k\delta)\right|\right) \le \Phi := \exp\left(T\left(3c+C\right)\rho_{\max}\right)\left(\|\rho_0''\|_\infty + \frac{\rho_{\max} W\left(\rho_{\max}+Y\right) + 3Y\gamma(\rho_{\max}) + 3CY^2}{\left(3c+C\right)\rho_{\max}}\right) \quad (5.43)$$

Equalities (5.39), (5.40), (5.41) in conjunction with (2.9), (2.11), (5.5), (5.13), (5.16), (5.21), (5.22), (5.23), (5.43) and the facts that $0 \le v_i(k\delta) \le v_{\max}$, $\left|v_{i+1}(k\delta) - v_i(k\delta)\right| \le hc\rho_{\max}$ for all $i = 0,\pm1,\pm2,...$, $k = 0,1,...,m$, give the following estimates:

$$\max_{i=...,-1,0,1,...}\left(\left|\psi_i(k\delta)\right|\right) \le \Psi := v_{\max}\Phi + \left(\gamma(\rho_{\max}) + (2c+C)Y\right)\rho_{\max}, \text{ for } k = 0,1,...,m-1 \quad (5.44)$$

$$\max_{i=...,-1,0,1,...}\left(\left|\eta_i(k\delta)\right|\right) \le H := \left(\gamma(\rho_{\max}) + (2c+C)Y\right)\rho_{\max} + v_{\max}\Phi, \text{ for } k = 0,1,...,m \quad (5.45)$$

$$\max_{i=...,-1,0,1,...}\left(\left|\mu_i(k\delta)\right|\right) \le M := \rho_{\max}\left(cG + (\lambda^{-1}+G)\Gamma(\rho_{\max})\right) + GLY + v_{\max}\Psi, \text{ for } k = 0,1,...,m-1 \quad (5.46)$$

Using (5.37), (5.19), (5.21), we obtain for $i = 0,\pm1,\pm2,...$, $k = 0,1,...,m$:

$$\left|\theta_{i+1}(k\delta) - \theta_i(k\delta)\right| \le h\left(\gamma(\rho_{\max}) + CY\right) \quad (5.47)$$

Using (5.37), (2.9), (5.5), (5.6), (5.16), (5.22), we obtain $i = 0,\pm1,\pm2,...$, $k = 0,1,...,m-1$:

$$\left|\theta_i((k+1)\delta) - \theta_i(k\delta)\right| \le \delta\Gamma(\rho_{\max})\left(\lambda^{-1} + G\right) \quad (5.48)$$

We define the functions $y(t,x;N)$, $g(t,x;N)$ for $(t,x) \in [0,T]\times\Re$ and for every integer $N > 2$ (recall that $h = N^{-1}$, $\delta = \lambda h$, $m\delta = T$):

$$\begin{gathered} y(k\delta,x;N) = \left(i+1-xN\right)y_i(k\delta) + \left(xN-i\right)y_{i+1}(k\delta), \\ g(k\delta,x;N) = \left(i+1-xN\right)g_i(k\delta) + \left(xN-i\right)g_{i+1}(k\delta), \\ \theta(k\delta,x;N) = \left(i+1-xN\right)\theta_i(k\delta) + \left(xN-i\right)\theta_{i+1}(k\delta) \\ \text{with } i = [xN], \text{ for } x \in \Re,\ k = 0,...,m, \end{gathered} \quad (5.49)$$

$$\begin{gathered} y(t,x;N) = \left(k+1-\lambda^{-1}tN\right)y(k\delta,x;N) + \left(\lambda^{-1}tN - k\right)y((k+1)\delta,x;N), \\ g(t,x;N) = \left(k+1-\lambda^{-1}tN\right)g(k\delta,x;N) + \left(\lambda^{-1}tN - k\right)g((k+1)\delta,x;N), \\ \theta(t,x;N) = \left(k+1-\lambda^{-1}tN\right)\theta(k\delta,x;N) + \left(\lambda^{-1}tN - k\right)\theta((k+1)\delta,x;N) \\ \text{with } k = [\lambda^{-1}tN] \text{ for } x \in \Re,\ t \in [0,T). \end{gathered} \quad (5.50)$$

It follows from (5.15), (5.16) and the fact that $\rho_{i+N}(k\delta) = \rho_i(k\delta)$, $v_{i+N}(k\delta) = v_i(k\delta)$ for all $i = 0,\pm1,\pm2,...$, $k = 0,1,...,m+1$, that $y_{i+N}(k\delta) = y_i(k\delta)$, $g_{i+N}(k\delta) = g_i(k\delta)$, $\theta_{i+N}(k\delta) = \theta_i(k\delta)$ for all $i = 0,\pm1,\pm2,...$, $k = 0,1,...,m$. Therefore, definitions (5.49), (5.50) imply that for each $t \in [0,T]$ the functions $y(t,\cdot;N)$, $g(t,\cdot;N)$, $\theta(t,\cdot;N)$ are periodic with period 1. Estimates (5.21), (5.22) and the fact that $\left|v_{i+1}(k\delta) - v_i(k\delta)\right| \le hc\rho_{\max}$ for all $i = 0,\pm1,\pm2,...$, $k = 0,1,...,m$, in conjunction with definitions (5.49), (5.50) imply that the following estimates hold for every integer $N > 2$:

$$\left|y(t,x;N)\right| \le Y \text{, for } (t,x)\in[0,T]\times\Re \tag{5.51}$$

$$\left|g(t,x;N)\right| \le G \text{, for } (t,x)\in[0,T]\times\Re \tag{5.52}$$

$$\left|\theta(t,x;N)\right| \le c\rho_{\max} \text{, for } (t,x)\in[0,T]\times\Re \tag{5.53}$$

Definitions (5.33), (5.37), (5.34), (5.35), (5.36) in conjunction with estimates (5.43), (5.44), (5.45), (5.46), (5.47), (5.48) imply the existence of a constant $\tilde{L}>0$ independent of $N>2$ for which the following estimate holds for all $i,j=0,\pm1,\pm2,...$, $k,l=0,1,...,m$:

$$\left|y_i(k\delta)-y_j(l\delta)\right|+\left|g_i(k\delta)-g_j(l\delta)\right|+\left|\theta_i(k\delta)-\theta_j(l\delta)\right|\le\tilde{L}(h\left|i-j\right|+\delta\left|k-l\right|) \tag{5.54}$$

Estimate (5.54) in conjunction with definitions (5.49), (5.50) implies that there exists a constant $L_2>0$ (independent of $N>2$) such that the following Lipschitz inequality holds for every integer $N>2$:

$$\left|y(t,x;N)-y(\tau,z;N)\right|+\left|g(t,x;N)-g(\tau,z;N)\right|+\left|\theta(t,x;N)-\theta(\tau,z;N)\right|\le L_2\left(\left|x-z\right|+\left|t-\tau\right|\right),$$
$$\text{for all } t,\tau\in[0,T] \text{ and } x,z\in\Re. \tag{5.55}$$

It follows from (5.51), (5.52), (5.55) that the sequences of functions $\{y(\cdot;N)\}_{N=3}^{\infty}$, $\{g(\cdot;N)\}_{N=3}^{\infty}$, $\{\theta(\cdot;N)\}_{N=3}^{\infty}$ are uniformly bounded and equicontinuous. Therefore, compactness of $[0,T]\times[0,1]$ and the Arzela-Ascoli theorem implies that there exist Lipschitz functions $y:[0,T]\times[0,1]\to\Re$, $g:[0,T]\times[0,1]\to\Re$, $\theta:[0,T]\times[0,1]\to\Re$ and subsequences $\left\{v(\cdot;N_q)\right\}_{q=1}^{\infty}$, $\left\{\rho(\cdot;N_q)\right\}_{q=1}^{\infty}$, $\left\{y(\cdot;N_q)\right\}_{q=1}^{\infty}$, $\left\{g(\cdot;N_q)\right\}_{q=1}^{\infty}$, $\left\{\theta(\cdot;N_q)\right\}_{q=1}^{\infty}$ for an increasing index sequence $\left\{N_q\right\}_{q=1}^{\infty}$, which converge uniformly on $[0,T]\times[0,1]$ to $v,\rho,y,g,\theta$, respectively. Moreover, the functions $y,g,\theta$ satisfy the same bounds with $y(\cdot;N),g(\cdot;N)$, i.e., $\left|y(t,x)\right|\le Y$, $\left|g(t,x)\right|\le G$, $\left|\theta(t,x)\right|\le c\rho_{\max}$, for $(t,x)\in[0,T]\times[0,1]$.

Since the functions $y(t,\cdot;N)$, $g(t,\cdot;N)$, $\theta(t,\cdot;N)$ are periodic with period 1, it follows that $y(t,1)=y(t,0)$, $\theta(t,1)=\theta(t,0)$ and $g(t,1)=g(t,0)$ for all $t\in[0,T]$. Therefore, the subsequences $\left\{y(\cdot;N_q)\right\}_{q=1}^{\infty}$, $\left\{g(\cdot;N_q)\right\}_{q=1}^{\infty}$, $\left\{\theta(\cdot;N_q)\right\}_{q=1}^{\infty}$ converge uniformly on $[0,T]\times\Re$ to the periodic extensions with respect to $x$ (with period 1) of $y,g,\theta$, respectively. We will denote by $y,g,\theta$ the periodic extensions with respect to $x$ (with period 1) of $y,g,\theta$ (again a slight abuse of notation).

We next show that $y(t,x)=\frac{\partial\rho}{\partial x}(t,x)$ for all $(t,x)\in[0,T]\times\Re$. Equivalently, we show that $\rho(t,x)-\rho(t,0)=\int_0^x y(t,z)dz$, for all $(t,x)\in[0,T]\times[0,1]$. Let $(t,x)\in[0,T]\times[0,1]$ be given (arbitrary). Using definitions (5.5), (5.24), (5.49), (5.15), inequalities (5.29), (5.55), (5.51), we obtain for $k=[\lambda^{-1}tN]$ and $i=[xN]$:

$$\left|\rho(t,x;N)-\rho(t,0;N)-\int_0^x y(t,z;N)dz\right| \le \left|\rho(t,x;N)-\rho(k\delta,x;N)\right|+\left|\rho(k\delta,x;N)-\rho(k\delta,ih;N)\right|$$
$$+\left|\rho(k\delta,0;N)-\rho(t,0;N)\right|+\left|\int_{ih}^x y(t,z;N)dz\right|+\left|\rho_i(k\delta)-\rho_0(k\delta)-\int_0^{ih} y(t,z;N)dz\right|$$
$$\le L_1(2t-2k\delta+x-ih)+Y(x-ih)+\left|h\sum_{s=0}^{i-1} y_s(k\delta)-\int_0^{ih} y(t,z;N)dz\right|$$
$$\le L_1(2\lambda+1)h+Yh+\left|\sum_{s=0}^{i-1}\int_{sh}^{(s+1)h} y(k\delta,sh;N)dz-\int_0^{ih} y(t,z;N)dz\right|$$
$$\le \left(L_1(2\lambda+1)+Y\right)h+\left|\sum_{s=0}^{i-1}\int_{sh}^{(s+1)h}\left(y(k\delta,sh;N)-y(t,z;N)\right)dz\right|$$
$$\le \left(L_1(2\lambda+1)+Y\right)h+\sum_{s=0}^{i-1}\int_{sh}^{(s+1)h}\left|y(k\delta,sh;N)-y(t,z;N)\right|dz$$
$$\le \left(L_1(2\lambda+1)+Y\right)h+\left(L_2(t-k\delta)+L_2h\right)ih \le \left(L_1(2\lambda+1)+Y+L_2\left(\lambda+1\right)\right)N^{-1}$$

In the above derivation, we have used the facts that $t-k\delta\le\delta$ and $x-ih\le h$. Since $\{y(\cdot;N_q)\}_{q=1}^{\infty}$, $\{\rho(\cdot;N_q)\}_{q=1}^{\infty}$ converge uniformly to $y$ and $q$ as $q\to+\infty$ (and $N_q\to+\infty$), the above inequality shows that $\rho(t,x)-\rho(t,0)=\int_0^x y(t,z)dz$ for all $(t,x)\in[0,T]\times[0,1]$.

A similar procedure shows that $g(t,x)=\frac{\partial\rho}{\partial t}(t,x)$ and $\theta(t,x)=\frac{\partial v}{\partial x}(t,x)$ for all $(t,x)\in[0,T]\times\Re$.

We next show that $g(t,x)+v(t,x)y(t,x)+\theta(t,x)\rho(t,x)=0$ for all $(t,x)\in[0,T]\times\Re$. Let $(t,x)\in[0,T]\times[0,1]$ be given (arbitrary). Since $g_i(k\delta)=-\theta_i(k\delta)\rho_i(k\delta)-v_i(k\delta)y_{i-1}(k\delta)$ for $i=0,\pm1,\pm2,...$, $k=0,1,...,m$ (a consequence of (5.18) and (5.37)), we obtain using (5.5), (5.29), (5.24), (5.26), (5.27), (5.49), (5.51) for $k=[\lambda^{-1}tN]$ and $i=\left[xN\right]$:

$$\left|g(t,x;N)+v(t,x;N)y(t,x;N)+\theta(t,x;N)\rho(t,x;N)\right|$$
$$\le\left|g(t,x;N)-g_i(k\delta)\right|+\left|v(t,x;N)y(t,x;N)-v_i(k\delta)y_{i-1}(k\delta)\right|+\left|\theta(t,x;N)\rho(t,x;N)-\theta_i(k\delta)\rho_i(k\delta)\right|$$
$$\le\left(L_2(1+\rho_{\max}+v_{\max})+L_1(Y+c\rho_{\max})\right)\left(t-k\delta+x-ih\right)+L_2v_{\max}h$$
$$\le\left(L_2(1+\rho_{\max}+v_{\max})+L_1(Y+c\rho_{\max})\right)\left(\delta+h\right)+L_2v_{\max}h$$
$$\le\left(\left(L_2(1+\rho_{\max}+v_{\max})+L_1(Y+c\rho_{\max})\right)\left(\lambda+1\right)+L_2v_{\max}\right)N^{-1}$$

In the above derivation, we have used the facts that $t-k\delta\le\delta$ and $x-ih\le h$. Since $\{v(\cdot;N_q)\}_{q=1}^{\infty}$, $\{\rho(\cdot;N_q)\}_{q=1}^{\infty}$, $\{y(\cdot;N_q)\}_{q=1}^{\infty}$, $\{g(\cdot;N_q)\}_{q=1}^{\infty}$, $\{\theta(\cdot;N_q)\}_{q=1}^{\infty}$ converge uniformly to $v,\rho,y,g,\theta$ as $q\to+\infty$ (and $N_q\to+\infty$), the above inequality shows that $g(t,x)+v(t,x)y(t,x)+\theta(t,x)\rho(t,x)=0$ for all $(t,x)\in[0,T]\times[0,1]$. Periodicity of $v,\rho,y,g,\theta$ implies that the equality $g(t,x)+v(t,x)y(t,x)+\theta(t,x)\rho(t,x)=0$ for all $(t,x)\in[0,T]\times\Re$.

Since $T > 0$ is arbitrary, we conclude that there exists a solution $\rho \in C^1(\Re_+ \times \Re)$ of the initial-value problem (2.1), (2.4), (2.5), (2.7) with $\rho[t] \in W^{2,\infty}(\Re) \cap Per(\Re)$ for all $t \geq 0$, which satisfies inequality (2.16) for all $t \geq 0, x \in \Re$.

Uniqueness follows by defining

$$E(t) := \frac{1}{2}\int_0^1 e^2(t,x)dx \tag{5.56}$$

where $e(t,x) := \rho(t,x) - \bar{\rho}(t,x)$ and $\rho, \bar{\rho} \in C^1(\Re_+ \times \Re)$ are solutions of the initial-value problem (2.1), (2.4), (2.5), (2.7) with $\rho[t], \bar{\rho}[t] \in W^{2,\infty}(\Re) \cap Per(\Re)$ for all $t \geq 0$. Let $T > 0$ be given (arbitrary) and define

$$\bar{\mu} := \max_{t\in[0,T]}\left(\left\|\frac{\partial v}{\partial x}[t]\right\|_\infty\right) + 2a\left(\max_{t\in[0,T]}\left(\|\rho[t]\|_\infty\right) + \max_{t\in[0,T]}\left(\|\bar{\rho}[t]\|_\infty\right)\right)\left(\max_{t\in[0,T]}\left(\left\|\frac{\partial \bar{\rho}}{\partial x}[t]\right\|_\infty\right) + \max_{t\in[0,T]}\left(\|\bar{\rho}[t]\|_\infty\right)\right) \tag{5.57}$$

where $a : \Re_+ \to \Re_+$ is the non-decreasing function involved in (2.6). Using (2.6) and (5.56), (5.57), we have for all $t \in [0,T]$:

$$\begin{aligned}
\dot{E}(t) &= -\int_0^1 e(t,x)\frac{\partial}{\partial x}\left(\rho(t,x)v(t,x) - \bar{\rho}(t,x)\bar{v}(t,x)\right)dx \\
&= -\frac{1}{2}\int_0^1 e^2(t,x)\frac{\partial v}{\partial x}(t,x)dx - \int_0^1 e(t,x)\left(v(t,x) - \bar{v}(t,x)\right)\frac{\partial \bar{\rho}}{\partial x}(t,x)dx \\
&\quad -\int_0^1 e(t,x)\bar{\rho}(t,x)\frac{\partial}{\partial x}\left(v(t,x) - \bar{v}(t,x)\right)dx \leq \bar{\mu}E(t)
\end{aligned} \tag{5.58}$$

The differential inequality (5.58) implies that $E(t) \leq \exp(\bar{\mu}t)E(0)$ for all $t \in [0,T]$. Since $E(0) = 0$, we get $E(t) = 0$ for all $t \in [0,T]$. Since $e[t]$ is periodic with period 1 for all $t \geq 0$ and since $T > 0$ is arbitrary, we conclude that $\rho(t,x) \equiv \bar{\rho}(t,x)$. The proof is complete. $\lhd$

**Proof of Theorem 2.4:** It suffices to show that (2.6) holds for the mapping

$$(K(\rho))(x) = f\left(\int_x^{x+\eta} \omega(s-x)\rho(s)ds\right)g\left(\int_{x-\zeta}^{x} \tilde{\omega}(x-s)\rho(s)ds\right), \text{ for } \rho \in Per(\Re),\ x \in \Re \tag{5.59}$$

and that the parameterized family

$$K_N(\rho) = f\left(\sum_{i=0}^{N-1}\rho_i \int_{ih}^{(i+1)h} \omega(s)ds\right)g\left(\sum_{i=1}^{N-1}\rho_i \int_{1-ih}^{1-(i-1)h} \tilde{\omega}(s)ds\right), \text{ for } \rho \in \Re_+^N \tag{5.60}$$

where $h = N^{-1}$, smoothly approximates the mapping $K \in C^0(Per(\Re); Per(\Re))$ defined by (5.59).

First we show the validity of (2.6) for certain non-decreasing function $a : \Re_+ \to \Re_+$. For every $\rho, \bar{\rho} \in Per(\Re)$ and $x \in \Re$, definition (5.59) implies the existence of $\theta, \vartheta \geq 0$ such that

$$\begin{aligned}
(K(\rho))(x) - (K(\bar{\rho}))(x) &= f'(\theta)\left(\int_x^{x+\eta} \omega(s-x)\left(\rho(s) - \bar{\rho}(s)\right)ds\right)g\left(\int_{x-\zeta}^{x} \tilde{\omega}(x-s)\rho(s)ds\right) \\
&+ f\left(\int_x^{x+\eta} \omega(s-x)\bar{\rho}(s)ds\right)g'(\vartheta)\left(\int_{x-\zeta}^{x} \tilde{\omega}(x-s)\left(\rho(s) - \bar{\rho}(s)\right)ds\right)
\end{aligned}$$

The above equation in conjunction with (2.17) and the facts that $\eta,\zeta\in(0,1]$ and that $\omega,\tilde{\omega}$ are non-increasing functions, allow us to obtain by using the Cauchy-Schwarz inequality, the following estimate for all $\rho,\bar{\rho}\in Per(\Re)$ and $x\in\Re$:

$$\left|(K(\rho))(x)-(K(\bar{\rho}))(x)\right|\leq M^2\omega(0)\int_x^{x+\eta}\left|\rho(s)-\bar{\rho}(s)\right|ds+M^2\tilde{\omega}(0)\int_{x-\zeta}^{x}\left|\rho(s)-\bar{\rho}(s)\right|ds$$

$$\leq M^2\left(\omega(0)+\tilde{\omega}(0)\right)\int_0^1\left|\rho(s)-\bar{\rho}(s)\right|ds\leq M^2\left(\omega(0)+\tilde{\omega}(0)\right)\left(\int_0^1\left|\rho(s)-\bar{\rho}(s)\right|^2 ds\right)^{1/2}$$

Using the above estimate we obtain for all $\rho,\bar{\rho}\in Per(\Re)$ and $x\in\Re$

$$\left|\rho(x)-\bar{\rho}(x)\right|\left|(K(\rho))(x)-(K(\bar{\rho}))(x)\right|\leq\frac{1}{2}\left|\rho(x)-\bar{\rho}(x)\right|^2+\frac{1}{2}M^4\left(\omega(0)+\tilde{\omega}(0)\right)^2\int_0^1\left|\rho(s)-\bar{\rho}(s)\right|^2 ds$$

which directly implies the following estimate for all $\rho,\bar{\rho}\in Per(\Re)$:

$$\int_0^1\left|\rho(x)-\bar{\rho}(x)\right|\left|(K(\rho))(x)-(K(\bar{\rho}))(x)\right|dx\leq\frac{1}{2}\left(1+M^4\left(\omega(0)+\tilde{\omega}(0)\right)^2\right)\int_0^1\left|\rho(s)-\bar{\rho}(s)\right|^2 ds \qquad (5.61)$$

Notice that definition (5.59) implies the following equation for all $\rho\in Per(\Re)$ and $x\in\Re$:

$$\frac{\partial}{\partial x}(K(\rho))(x)=$$

$$\left(\omega(\eta)\rho(x+\eta)-\omega(0)\rho(x)-\int_x^{x+\eta}\omega'(s-x)\rho(s)ds\right)f'\left(\int_x^{x+\eta}\omega(s-x)\rho(s)ds\right)g\left(\int_{x-\zeta}^{x}\tilde{\omega}(x-s)\rho(s)ds\right) \qquad (5.62)$$

$$+\left(\tilde{\omega}(0)\rho(x)-\tilde{\omega}(\zeta)\rho(x-\zeta)+\int_{x-\zeta}^{x}\tilde{\omega}'(x-s)\rho(s)ds\right)f\left(\int_x^{x+\eta}\omega(s-x)\rho(s)ds\right)g'\left(\int_{x-\zeta}^{x}\tilde{\omega}(x-s)\rho(s)ds\right)$$

The Cauchy-Schwarz inequality and the fact that $\eta\in(0,1]$ gives the following estimate for all $\rho,\bar{\rho}\in Per(\Re)$ and $x\in\Re$

$$A(x):=\left|\omega(\eta)\left(\rho(x+\eta)-\bar{\rho}(x+\eta)\right)-\omega(0)\left(\rho(x)-\bar{\rho}(x)\right)-\int_x^{x+\eta}\omega'(s-x)\left(\rho(s)-\bar{\rho}(s)\right)ds\right|$$

$$\leq\omega(\eta)\left|\rho(x+\eta)-\bar{\rho}(x+\eta)\right|+\omega(0)\left|\rho(x)-\bar{\rho}(x)\right|+\int_x^{x+\eta}\left|\omega'(s-x)\right|\left|\rho(s)-\bar{\rho}(s)\right|ds$$

$$\leq\omega(\eta)\left|\rho(x+\eta)-\bar{\rho}(x+\eta)\right|+\omega(0)\left|\rho(x)-\bar{\rho}(x)\right|+\left(\int_0^{\eta}\left|\omega'(s)\right|^2 ds\right)^{1/2}\left(\int_x^{x+\eta}\left|\rho(s)-\bar{\rho}(s)\right|^2 ds\right)^{1/2}$$

$$\leq\omega(\eta)\left|\rho(x+\eta)-\bar{\rho}(x+\eta)\right|+\omega(0)\left|\rho(x)-\bar{\rho}(x)\right|+\left(\int_0^{\eta}\left|\omega'(s)\right|^2 ds\right)^{1/2}\left(\int_0^{1}\left|\rho(s)-\bar{\rho}(s)\right|^2 ds\right)^{1/2}$$

which implies that

$$A(x)\left|\rho(x)-\bar{\rho}(x)\right|\leq\frac{1}{2}\omega^2(\eta)\left|\rho(x+\eta)-\bar{\rho}(x+\eta)\right|^2$$

$$+(1+\omega(0))\left|\rho(x)-\bar{\rho}(x)\right|^2+\frac{1}{2}\left(\int_0^{\eta}\left|\omega'(s)\right|^2 ds\right)\left(\int_0^{1}\left|\rho(s)-\bar{\rho}(s)\right|^2 ds\right)$$

and consequently

$$\int_0^1 A(x)\left|\rho(x)-\bar{\rho}(x)\right|dx \le \left(\frac{1}{2}\omega^2(\eta)+1+\omega(0)+\frac{1}{2}\int_0^\eta |\omega'(s)|^2\,ds\right)\left(\int_0^1 |\rho(s)-\bar{\rho}(s)|^2\,ds\right)$$

Obtaining a similar estimate for the term

$$B(x):=\left|\tilde{\omega}(0)\big(\rho(x)-\bar{\rho}(x)\big)-\tilde{\omega}(\zeta)\big(\rho(x-\zeta)-\bar{\rho}(x-\zeta)\big)+\int_{x-\zeta}^{x}\tilde{\omega}'(x-s)\big(\rho(s)-\bar{\rho}(s)\big)ds\right|$$

and noticing that there exists a constant $\Theta>0$ such that $\left|\omega(\eta)\rho(x+\eta)-\omega(0)\rho(x)-\int_x^{x+\eta}\omega'(s-x)\rho(s)ds\right|\le\Theta\|\rho\|_\infty$, $\left|\tilde{\omega}(0)\rho(x)-\tilde{\omega}(\zeta)\rho(x-\zeta)+\int_{x-\zeta}^{x}\tilde{\omega}'(x-s)\rho(s)ds\right|\le\Theta\|\rho\|_\infty$ for all $x\in\Re$ and $\rho\in Per(\Re)$, we conclude from (5.62), by using the same arguments as above (for the derivation of (5.61)), that there exists a constant $\bar{\Theta}>0$ which satisfies the following estimate for all $\rho,\bar{\rho}\in Per(\Re)$:

$$\int_0^1 |\rho(x)-\bar{\rho}(x)|\left|\frac{\partial}{\partial x}\big((K(\rho))(x)-(K(\bar{\rho}))(x)\big)\right|dx \le \bar{\Theta}\big(1+\|\rho\|_\infty+\|\bar{\rho}\|_\infty\big)\int_0^1|\rho(x)-\bar{\rho}(x)|^2\,dx \tag{5.63}$$

Inequality (2.6) with an appropriate (linear) non-decreasing function $a:\Re_+\to\Re_+$ is a direct consequence of estimates (5.61) and (5.63).

By virtue of Lemma 2.2, in order to show that the parameterized family $K_N$ defined by (5.60) smoothly approximates the mapping $K\in C^0(Per(\Re);Per(\Re))$ defined by (5.59), it suffices to show that the parameterized families

$$F_N(\rho)=f\left(B_N\rho\right),\ \text{for}\ \rho\in\Re_+^N \tag{5.64}$$

$$G_N(\rho)=g\left(\bar{B}_N\rho\right),\ \text{for}\ \rho\in\Re_+^N \tag{5.65}$$

where

$$B_N\rho=\sum_{i=0}^{N-1}\rho_i\int_{ih}^{(i+1)h}\omega(s)ds\ ,\ \text{for}\ \rho\in\Re_+^N \tag{5.66}$$

$$\bar{B}_N\rho=\sum_{i=1}^{N-1}\rho_i\int_{1-ih}^{1-(i-1)h}\tilde{\omega}(s)ds\ ,\ \text{for}\ \rho\in\Re_+^N \tag{5.67}$$

$h=N^{-1}$, smoothly approximate the mappings

$$(F(\rho))(x)=f\left(\int_x^{x+1}\omega(s-x)\rho(s)ds\right),\ \text{for}\ \rho\in Per(\Re),\ x\in\Re \tag{5.68}$$

$$(G(\rho))(x)=g\left(\int_{x-1}^{x}\tilde{\omega}(x-s)\rho(s)ds\right),\ \text{for}\ \rho\in Per(\Re),\ x\in\Re \tag{5.69}$$

Using the facts that $h=1/N$, $\rho^{(i)}=(\rho_i,...,\rho_{N-1},\rho_0,...,\rho_{i-1})$ for $i=1,...,N-1$, $y=(y_0,...,y_{N-1})^\top=h^{-1}(\rho^{(1)}-\rho)$ and $\varphi=(\varphi_0,...,\varphi_{N-1})^\top=h^{-2}(\rho^{(2)}-2\rho^{(1)}+\rho)$, we establish the following useful identities hold for the linear mappings $B_N\rho$, $\bar{B}_N\rho$ defined by (5.66), (5.67) for all $\rho\in\Re_+^N$:

$$B_N\rho^{(1)}-B_N\rho=\sum_{i=1}^{N-1}\rho_i\int_{(i-1)h}^{ih}\big(\omega(s)-\omega(s+h)\big)ds-\rho_0\left(\int_0^h\omega(s)ds-\int_{1-h}^{1}\omega(s)ds\right) \tag{5.70}$$

$$\bar{B}_N\rho^{(1)}-\bar{B}_N\rho=\rho_0\int_h^{2h}\tilde{\omega}(s)ds-\sum_{i=2}^{N-1}\rho_i\int_{1-ih}^{1-(i-1)h}\left(\tilde{\omega}(s)-\tilde{\omega}(s+h)\right)ds-\rho_1\int_{1-h}^{1}\tilde{\omega}(s)ds \quad (5.71)$$

$$2B_N\rho^{(1)}-B_N\rho-B_N\rho^{(2)}=hy_0\left(\int_0^h\omega(s)ds-\int_{1-h}^{1}\omega(s)ds\right)-h\sum_{i=1}^{N-1}y_i\int_{(i-1)h}^{ih}\left(\omega(s)-\omega(s+h)\right)ds \quad (5.72)$$

$$2\bar{B}_N\rho^{(1)}-\bar{B}_N\rho-\bar{B}_N\rho^{(2)}=h\sum_{i=2}^{N-1}y_i\int_{1-ih}^{1-(i-1)h}\left(\tilde{\omega}(s)-\tilde{\omega}(s+h)\right)ds+hy_1\int_{1-h}^{1}\tilde{\omega}(s)ds-hy_0\int_h^{2h}\tilde{\omega}(s)ds \quad (5.73)$$

$$\begin{aligned}&3B_N\rho^{(1)}+B_N\rho^{(3)}-B_N\rho-3B_N\rho^{(2)}\\&=h^2\sum_{i=1}^{N-1}\varphi_i\int_{(i-1)h}^{ih}\left(\omega(s)-\omega(s+h)\right)ds-h^2\varphi_0\left(\int_0^h\omega(s)ds-\int_{1-h}^{1}\omega(s)ds\right)\end{aligned} \quad (5.74)$$

$$\begin{aligned}&3\bar{B}_N\rho^{(1)}+\bar{B}_N\rho^{(3)}-\bar{B}_N\rho-3\bar{B}_N\rho^{(2)}=\\&h^2\varphi_0\int_h^{2h}\tilde{\omega}(s)ds-h^2\varphi_1\int_{1-h}^{1}\tilde{\omega}(s)ds-h^2\sum_{i=2}^{N-1}\varphi_i\int_{1-ih}^{1-(i-1)h}\left(\tilde{\omega}(s)-\tilde{\omega}(s+h)\right)ds\end{aligned} \quad (5.75)$$

Notice that since $\omega,\tilde{\omega}$ are non-negative and non-increasing functions, it follows from identities (2.10), (5.71) that the following estimates hold for every $0<\rho_{\min}\le\rho_{\max}$ and $\rho\in\Re_+^N$ with $\rho_{\min}1_N\le\rho\le\rho_{\max}1_N$ :

$$\left(\int_0^h\omega(s)ds-\int_{1-h}^{1}\omega(s)ds\right)\left(\rho_{\min}-\rho_0\right)\le B_N\rho^{(1)}-B_N\rho\le\left(\int_0^h\omega(s)ds-\int_{1-h}^{1}\omega(s)ds\right)\left(\rho_{\max}-\rho_0\right) \quad (5.76)$$

$$\left(\rho_0-\rho_{\max}\right)\int_h^{2h}\tilde{\omega}(s)ds\le\bar{B}_N\rho^{(1)}-\bar{B}_N\rho\le\left(\rho_0-\rho_{\min}\right)\int_h^{2h}\tilde{\omega}(s)ds \quad (5.77)$$

$$\left|B_N\rho^{(1)}-B_N\rho\right|\le h\omega(0)\rho_{\max} \quad (5.78)$$

$$\left|\bar{B}_N\rho^{(1)}-\bar{B}_N\rho\right|\le h\tilde{\omega}(0)\rho_{\max} \quad (5.79)$$

Inequalities (2.10) for $F_N(\rho)$ and $G_N(\rho)$ (with appropriate $c$) are consequences of estimates (5.76), (5.77), definitions (5.64), (5.65), the fact that $f'(\rho)\le 0$, $g'(\rho)\ge 0$ for all $\rho\ge 0$, inequality (2.17) and the mean value theorem for $f$ and $g$.

Inequality (2.17) implies the following estimates for all $x_1,x_2,z_1,z_2\in\Re_+$:

$$\begin{aligned}&\left|f(x_2)-f(x_1)-f(z_2)+f(z_1)\right|\le M\left|x_2-x_1-z_2+z_1\right|+M\left|x_2-x_1\right|\left|x_1-z_1\right|\\&+M\left|x_2-x_1\right|\left|z_2-z_1\right|+M\left|x_2-x_1\right|^2\end{aligned} \quad (5.80)$$

$$\begin{aligned}&\left|g(x_2)-g(x_1)-g(z_2)+g(z_1)\right|\le M\left|x_2-x_1-z_2+z_1\right|+M\left|x_2-x_1\right|\left|x_1-z_1\right|\\&+M\left|x_2-x_1\right|\left|z_2-z_1\right|+M\left|x_2-x_1\right|^2\end{aligned} \quad (5.81)$$

Indeed, by virtue of the mean value theorem for all $x_1,x_2,z_1,z_2\in\Re_+$ there exist $\lambda,\mu\in(0,1)$ such that

$$\begin{aligned} &f(x_2)-f(x_1)-f(z_2)+f(z_1)=f'(\theta)(x_2-x_1)-f'(\vartheta)(z_2-z_1) \\ &=f'(\vartheta)(x_2-x_1-z_2+z_1)+f''(r)(\theta-\vartheta)(x_2-x_1) \end{aligned}$$

where $\theta = x_1+\lambda(x_2-x_1)$, $\vartheta = z_1+\mu(z_2-z_1)$ and $r\in\Re$. A similar relation holds for $g$ as well. Inequalities (5.80), (5.81) are consequences of the above equality and estimate (2.17).

Let $0<\rho_{\min}\le\rho_{\max}$ be given and let arbitrary vectors $\rho,\tilde{\rho}\in\Re_+^N$ with $\rho_{\min}1_N\le\rho\le\rho_{\max}1_N$, $\rho_{\min}1_N\le\tilde{\rho}\le\rho_{\max}1_N$ be also given. Using (5.80), (5.81) with $x_2=B_N\rho$, $x_1=B_N\rho^{(1)}$, $z_2=B_N\tilde{\rho}$, $z_1=B_N\tilde{\rho}^{(1)}$ for $f$ and $x_2=\bar{B}_N\rho$, $x_1=\bar{B}_N\rho^{(1)}$, $z_2=\bar{B}_N\tilde{\rho}$, $z_1=\bar{B}_N\tilde{\rho}^{(1)}$ for $g$ in conjunction with (5.78), (5.79), (5.64), (5.65) and the following inequalities

$$\left|B_N\rho-B_N\tilde{\rho}\right|\le\left|\rho-\tilde{\rho}\right|_\infty\int_0^1\omega(s)ds \tag{5.82}$$

$$\left|\bar{B}_N\rho-\bar{B}_N\tilde{\rho}\right|\le\left|\rho-\tilde{\rho}\right|_\infty\int_0^1\tilde{\omega}(s)ds \tag{5.83}$$

$$\left|B_N\rho^{(1)}-B_N\rho-B_N\tilde{\rho}^{(1)}+B_N\tilde{\rho}\right|\le 2\left|\rho-\tilde{\rho}\right|_\infty\left(\int_0^h\omega(s)ds-\int_{1-h}^1\omega(s)ds\right)$$

$$\left|\bar{B}_N\rho^{(1)}-\bar{B}_N\rho-\bar{B}_N\tilde{\rho}^{(1)}+\bar{B}_N\tilde{\rho}\right|\le 2\left|\rho-\tilde{\rho}\right|_\infty\int_h^{2h}\tilde{\omega}(s)ds$$

which are direct consequences of (5.66), (5.67), (5.70), (5.71) as well as the fact that $\omega,\tilde{\omega}$ are non-negative and non-increasing functions, we obtain:

$$\begin{aligned} &\left|F_N(\rho)-F_N(\rho^{(1)})-F_N(\tilde{\rho})+F_N(\tilde{\rho}^{(1)})\right|\le \\ &Mh\omega(0)\left|\rho-\tilde{\rho}\right|_\infty\left(2+\rho_{\max}\int_0^1\omega(s)ds\right)+2Mh^2\omega^2(0)\rho_{\max}^2 \end{aligned}$$

$$\begin{aligned} &\left|G_N(\rho)-G_N(\rho^{(1)})-G_N(\tilde{\rho})+G_N(\tilde{\rho}^{(1)})\right|\le \\ &Mh\tilde{\omega}(0)\left|\rho-\tilde{\rho}\right|_\infty\left(2+\rho_{\max}\int_0^1\tilde{\omega}(s)ds\right)+2Mh^2\tilde{\omega}^2(0)\rho_{\max}^2 \end{aligned}$$

The above estimates directly imply the validity of (2.9) for $F_N(\rho)$, $G_N(\rho)$ and for an appropriate non-decreasing function $\Gamma:\Re_+\to\Re_+$. Moreover, inequalities (5.82), (5.83) in conjunction with the mean value theorem and (2.17) imply the validity of (2.8) for $F_N(\rho)$, $G_N(\rho)$ and for an appropriate constant $L>0$.

Inequality (2.17) implies that for every $b\ge a\ge 0$ and $x_1,x_2,x_3,x_4\in[a,b]$, the following inequalities hold:

$$\left|2f(x_2)-f(x_1)-f(x_3)\right|\le M\left|2x_2-x_1-x_3\right|+M(b-a)^2 \tag{5.84}$$

$$\left|2g(x_2)-g(x_1)-g(x_3)\right|\le M\left|2x_2-x_1-x_3\right|+M(b-a)^2 \tag{5.85}$$

$$\begin{aligned} &\left|3f(x_2)-3f(x_3)+f(x_4)-f(x_1)\right|\le M\left|3x_2+x_4-3x_3-x_1\right|+\frac{7}{6}M(b-a)^3 \\ &+3M(b-a)\left|2x_3-x_2-x_4\right|+M\left|2x_2-x_1-x_3\right|^2+M\left|2x_3-x_2-x_4\right|^2 \end{aligned} \tag{5.86}$$

$$\left|3g(x_2)-3g(x_3)+g(x_4)-g(x_1)\right| \le M\left|3x_2+x_4-3x_3-x_1\right|+\frac{7}{6}M(b-a)^3$$
$$+3M(b-a)\left|2x_3-x_2-x_4\right|+M\left|2x_2-x_1-x_3\right|^2+M\left|2x_3-x_2-x_4\right|^2 \tag{5.87}$$

Indeed, Taylor's theorem implies the existence of $\theta,\vartheta,r\in[a,b]$ so that

$$3f(x_2)-3f(x_3)+f(x_4)-f(x_1)=f'(x_1)(3x_2+x_4-3x_3-x_1)$$
$$+\frac{1}{2}f''(x_1)\left(3(x_2-x_1)^2-3(x_3-x_1)^2+(x_4-x_1)^2\right)$$
$$+\frac{1}{2}f'''(\theta)(x_2-x_1)^3-\frac{1}{2}f'''(\vartheta)(x_3-x_1)^3+\frac{1}{6}f'''(r)(x_4-x_1)^3$$

Using the above equality in conjunction with (2.17), the fact that $3(x_2-x_1)^2-3(x_3-x_1)^2+(x_4-x_1)^2=-6\kappa(x_3-x_2)+(s-\kappa)^2$, where $s=2x_2-x_1-x_3$, $\kappa=2x_3-x_2-x_4$ and the fact that $x_1,x_2,x_3,x_4\in[a,b]$, we obtain (5.86). A similar derivation gives inequalities (5.84), (5.85), (5.87).

Inequalities (5.72), (5.73), (5.84), (5.85) allow us to conclude the validity of (2.11) for $F_N(\rho)$, $G_N(\rho)$ with $C=2M\left(\tilde{\omega}(0)+\omega(0)\right)$ and for an appropriate non-decreasing function $\gamma:\Re_+\to\Re_+$. Indeed, using (5.72) and the fact that $\omega$ is a non-negative and non-increasing function, we get:

$$\left|2B_N\rho^{(1)}-B_N\rho-B_N\rho^{(2)}\right|\le 2h^2\omega(0)\left|y\right|_\infty \tag{5.88}$$

Applying (5.84) with $x_2=B_N\rho^{(1)}$, $x_1=B_N\rho$, $x_3=B_N\rho^{(2)}$ and definition (5.64) in conjunction with the above inequality, we get:

$$\left|2F_N(\rho^{(1)})-F_N(\rho)-F_N(\rho^{(2)})\right|\le 2Mh^2\omega(0)\left|y\right|_\infty+M(b-a)^2 \tag{5.89}$$

where $b\ge a\ge 0$ are any numbers for which $B_N\rho,B_N\rho^{(1)},B_N\rho^{(2)}\in[a,b]$. Since (5.78) implies that $\left|B_N\rho^{(1)}-B_N\rho\right|\le h\omega(0)\rho_{\max}$, $\left|B_N\rho^{(2)}-B_N\rho^{(1)}\right|\le h\omega(0)\rho_{\max}$, it follows that there exist numbers $b\ge a\ge 0$ with $b-a\le 2h\omega(0)\rho_{\max}$ such that $B_N\rho,B_N\rho^{(1)},B_N\rho^{(2)}\in[a,b]$. It follows from (5.89) that

$$\left|2F_N(\rho^{(1)})-F_N(\rho)-F_N(\rho^{(2)})\right|\le 2Mh^2\omega(0)\left(\left|y\right|_\infty+2\omega(0)\rho_{\max}^2\right) \tag{5.90}$$

A similar procedure allows us to show that

$$\left|2G_N(\rho^{(1)})-G_N(\rho)-G_N(\rho^{(2)})\right|\le 2Mh^2\tilde{\omega}(0)\left(\left|y\right|_\infty+2\tilde{\omega}(0)\rho_{\max}^2\right) \tag{5.91}$$

Inequalities (5.90), (5.91) show the validity of (2.11) for $F_N(\rho)$, $G_N(\rho)$ with $C=2M\left(\tilde{\omega}(0)+\omega(0)\right)$ and for an appropriate non-decreasing function $\gamma:\Re_+\to\Re_+$.

Inequalities (5.74), (5.75), (5.86), (5.87) allow us to conclude the validity of (2.12) for $F_N(\rho)$, $G_N(\rho)$ with $C=2M\left(\tilde{\omega}(0)+\omega(0)\right)$ and for an appropriate non-decreasing function $W:\Re_+\to\Re_+$. Indeed, using (5.74) and the fact that $\omega$ is a non-negative and non-decreasing function, we get:

$$\left|3B_N\rho^{(1)}+B_N\rho^{(3)}-B_N\rho-3B_N\rho^{(2)}\right|\le 2h^2\left|\varphi\right|_\infty\left(\int_0^h\omega(s)ds-\int_{1-h}^1\omega(s)ds\right)$$

Applying (5.84) with $x_2 = B_N \rho^{(1)}$, $x_1 = B_N \rho$, $x_3 = B_N \rho^{(2)}$, $x_4 = B_N \rho^{(3)}$ and definition (5.64) in conjunction with the above inequality, we get:

$$\begin{aligned}
&\left|3F_N(\rho^{(1)}) + F_N(\rho^{(3)}) - F_N(\rho) - 3F_N(\rho^{(2)})\right| \le 2Mh^3\omega(0)|\varphi|_\infty + \frac{7}{6}M(b-a)^3 \\
&+3M(b-a)\left|2B_N\rho^{(2)} - B_N\rho^{(1)} - B_N\rho^{(3)}\right| + M\left|2B_N\rho^{(1)} - B_N\rho - B_N\rho^{(2)}\right|^2 \\
&+M\left|2B_N\rho^{(2)} - B_N\rho^{(1)} - B_N\rho^{(3)}\right|^2
\end{aligned} \tag{5.92}$$

where $b \ge a \ge 0$ are any numbers for which $B_N\rho, B_N\rho^{(1)}, B_N\rho^{(2)}, B_N\rho^{(3)} \in [a,b]$. Using (5.88) and its direct consequence $\left|2B_N\rho^{(2)} - B_N\rho^{(1)} - B_N\rho^{(3)}\right| \le 2h^2\omega(0)\left|y^{(1)}\right|_\infty = 2h^2\omega(0)|y|_\infty$ in conjunction with (5.92), we get:

$$\begin{aligned}
&\left|3F_N(\rho^{(1)}) + F_N(\rho^{(3)}) - F_N(\rho) - 3F_N(\rho^{(2)})\right| \le 2Mh^3\omega(0)|\varphi|_\infty \\
&+\frac{7}{6}M(b-a)^3 + 6M(b-a)h^2\omega(0)|y|_\infty + 8Mh^4\omega^2(0)\left(|y|_\infty\right)^2
\end{aligned} \tag{5.93}$$

Since (5.78) implies that $\left|B_N\rho^{(1)} - B_N\rho\right| \le h\omega(0)\rho_{\max}$, $\left|B_N\rho^{(2)} - B_N\rho^{(1)}\right| \le h\omega(0)\rho_{\max}$, $\left|B_N\rho^{(3)} - B_N\rho^{(2)}\right| \le h\omega(0)\rho_{\max}$, it follows that there exist numbers $b \ge a \ge 0$ with $b - a \le 3h\omega(0)\rho_{\max}$ such that $B_N\rho, B_N\rho^{(1)}, B_N\rho^{(2)}, B_N\rho^{(3)} \in [a,b]$. It follows from (5.93) and the fact that $h \le 1$ that

$$\begin{aligned}
&\left|3F_N(\rho^{(1)}) + F_N(\rho^{(3)}) - F_N(\rho) - 3F_N(\rho^{(2)})\right| \le \\
&2Mh^3\omega(0)\left(|\varphi|_\infty + \frac{63}{4}\omega^2(0)\rho_{\max}^3 + 9\omega(0)\rho_{\max}|y|_\infty + 4\omega(0)\left(|y|_\infty\right)^2\right)
\end{aligned} \tag{5.94}$$

A similar procedure allows us to show that

$$\begin{aligned}
&\left|3G_N(\rho^{(1)}) + G_N(\rho^{(3)}) - G_N(\rho) - 3G_N(\rho^{(2)})\right| \le \\
&2Mh^3\tilde{\omega}(0)\left(|\varphi|_\infty + \frac{63}{4}\tilde{\omega}^2(0)\rho_{\max}^3 + 9\tilde{\omega}(0)\rho_{\max}|y|_\infty + 4\tilde{\omega}(0)\left(|y|_\infty\right)^2\right)
\end{aligned} \tag{5.95}$$

Inequalities (5.94), (5.95) show the validity of (2.12) for $F_N(\rho)$, $G_N(\rho)$ with $C = 2M\left(\tilde{\omega}(0) + \omega(0)\right)$ and for an appropriate non-decreasing function $W : \Re_+ \to \Re_+$.

Finally, we show the validity of (2.13) for $F_N(\rho)$, $G_N(\rho)$ for an appropriate constant $S > 0$. Definitions (2.14), (2.15) imply that

$$(P_N\rho)(x) = h^{-1}\sum_{i=0}^{N-1}\rho_i\phi(x - ih) \text{ for } x \in \Re \tag{5.96}$$

where $\phi : \Re \to \Re$ is the periodic extension (with period 1) of the function

$$\phi(x) = \begin{cases} h - x, & x \in [0,h] \\ x - 1 + h, & x \in [1-h,1] \\ 0, x \notin [0,h] \cup [1-h,1] \end{cases}, \text{ for } x \in [0,1] \tag{5.97}$$

Formula (5.96) and periodicity of $\phi : \Re \to \Re$ in conjunction with definitions (5.68), (5.69) imply that

$$(F(P_N\rho))(ih) = (F(P_N\rho^{(i)}))(0),\ (G(P_N\rho))(ih) = (G(P_N\rho^{(i)}))(0), \text{ for } i = 0,...,N-1 \tag{5.98}$$

Therefore, it suffices to show that there exists a constant $S>0$ so that (2.13) holds for $F_N(\rho)$ and $G_N(\rho)$ with $i=0$. Definitions (5.64), (5.65), (5.66), (5.67), (5.68), (5.69) in conjunction with (5.96) and inequality (2.17) imply the inequalities:

$$\left|(F(P_N\rho))(0)-F_N(\rho)\right|\leq Mh^{-1}\left|\sum_{j=0}^{N-1}\rho_j\int_0^1\omega(s)\phi(s-jh)ds-\sum_{j=0}^{N-1}\rho_j h\int_{jh}^{(j+1)h}\omega(s)ds\right| \tag{5.99}$$

$$\left|(G(\rho))(0)-G_N(\rho)\right|\leq h^{-1}M\left|\sum_{j=0}^{N-1}\rho_j\int_0^1\tilde{\omega}(s)\phi(-s-jh)ds-\sum_{j=1}^{N-1}\rho_j h\int_{i-jh}^{1-(j-1)h}\tilde{\omega}(s)ds\right| \tag{5.100}$$

Using (5.97) and the fact that $\phi:\Re\to\Re$ is periodic with period 1, we get from (5.99) and (5.100):

$$\begin{aligned}&\left|(F(P_N\rho))(0)-F_N(\rho)\right|\\&\leq Mh^{-1}\left|\sum_{j=1}^{N-1}\rho_j\int_{jh}^{(j+1)h}\left(\omega(s)-\omega(s-h)\right)(jh-s)ds+\rho_0\int_{1-h}^{1}\omega(s)(s-1+h)ds-\rho_0\int_0^h s\omega(s)ds\right|\end{aligned} \tag{5.101}$$

$$\begin{aligned}&\left|(G(\rho))(0)-G_N(\rho)\right|\\&\leq h^{-1}M\left|\sum_{j=1}^{N-1}\rho_j\int_{1-jh}^{1-(j-1)h}\left(\tilde{\omega}(s-h)-\tilde{\omega}(s)\right)(s+jh-1)ds+\rho_0\int_{1-h}^{1}\tilde{\omega}(s)(h-1+s)ds+\rho_0\int_0^h\tilde{\omega}(s)(h-s)ds\right|\end{aligned} \tag{5.102}$$

Since $\omega,\tilde{\omega}$ are non-negative and non-increasing functions, it follows from inequalities (5.101), (5.102) that the following estimates hold for every $0<\rho_{\min}\leq\rho_{\max}$ and $\rho\in\Re_+^N$ with $\rho_{\min}1_N\leq\rho\leq\rho_{\max}1_N$:

$$\begin{aligned}&\left|(F(P_N\rho))(0)-F_N(\rho)\right|\\&\leq Mh^{-1}\rho_{\max}\left(\sum_{j=1}^{N-1}\int_{jh}^{(j+1)h}\left(\omega(s-h)-\omega(s)\right)(s-jh)ds+\int_{1-h}^{1}\omega(s)(s-1+h)ds+\int_0^h s\omega(s)ds\right)\end{aligned} \tag{5.103}$$

$$\begin{aligned}&\left|(G(\rho))(0)-G_N(\rho)\right|\\&\leq h^{-1}M\rho_{\max}\left(\sum_{j=1}^{N-1}\int_{1-jh}^{1-(j-1)h}\left(\tilde{\omega}(s-h)-\tilde{\omega}(s)\right)(s+jh-1)ds+\int_{1-h}^{1}\tilde{\omega}(s)(h-1+s)ds+\int_0^h\tilde{\omega}(s)(h-s)ds\right)\end{aligned} \tag{5.104}$$

Using the fact that $\int_I p(x)q(x)dx\leq\sup_{x\in I}\left(q(x)\right)\int_I p(x)dx$, for every pair of piecewise continuous, non-negative functions $p,q:I\to\Re$ and for every interval $I\subset\Re$, we obtain from (5.103), (5.104):

$$\left|(F(P_N\rho))(0)-F_N(\rho)\right|\leq M\rho_{\max}\left(\sum_{j=1}^{N-1}\int_{jh}^{(j+1)h}\left(\omega(s-h)-\omega(s)\right)ds+\int_{1-h}^{1}\omega(s)ds+\int_0^h\omega(s)ds\right) \tag{5.105}$$

$$\left|(G(\rho))(0)-G_N(\rho)\right|\leq M\rho_{\max}\left(\sum_{j=1}^{N-1}\int_{1-jh}^{1-(j-1)h}\left(\tilde{\omega}(s-h)-\tilde{\omega}(s)\right)ds+\int_{1-h}^{1}\tilde{\omega}(s)ds+\int_0^h\tilde{\omega}(s)ds\right) \tag{5.106}$$

Evaluating the integrals in the right hand sides of (5.105), (5.106) and using the fact that $\omega,\tilde{\omega}$ are non-negative and non-increasing functions, we get:

$$\left|(F(P_N\rho))(0)-F_N(\rho)\right|\leq 2Mh\rho_{\max}\omega(0) \tag{5.107}$$

$$\left|(G(\rho))(0)-G_N(\rho)\right|\leq 2Mh\rho_{\max}\tilde{\omega}(0) \tag{5.108}$$

Inequalities (5.107), (5.108) show the validity of (2.12) for $F_N(\rho)$, $G_N(\rho)$ with $S=2M\left(\tilde{\omega}(0)+\omega(0)\right)$. The proof is complete. $\triangleleft$

**Proof of Theorem 3.2:** Let $\rho_0 \in W^{2,\infty}(\Re)\cap Per(\Re)$ be a given function that satisfies condition (3.2). Let $\rho\in C^1(\Re_+\times\Re)$ with $\rho[t]\in W^{2,\infty}(\Re)\cap Per(\Re)$ for all $t\geq 0$ be the unique solution of (2.1), (2.3), (2.4), (2.7). We define for $t\geq 0$:

$$V(t)=\int_0^1\left(\rho(t,x)\ln\left(\frac{\rho(t,x)}{\rho^*}\right)+\rho^*-\rho(t,x)\right)dx \tag{5.109}$$

The time derivative of $V(t)$ can be computed using (2.1), (2.3), (2.4), (5.109) and the facts that $\rho[t]$ is periodic with period 1, $\rho^*=\int_x^{x+1}\rho(t,s)ds$ for all $t\geq 0$, $x\in\Re$, $\zeta=1$, $\omega(x)=\eta^{-1}$ for $x\in[0,\eta]$ and $\omega(x)=0$ for $x>\eta$, $\tilde{\omega}(x)=1-x$ for $x\in[0,1]$ and $f'(\rho)\leq 0$ for all $\rho\geq 0$:

$$\begin{aligned}
\dot{V}(t)&=\int_0^1\left(\rho^*-\rho(t,x)\right)\frac{\partial v}{\partial x}(t,x)dx\\
&=\eta^{-1}\int_0^1\left(\rho^*-\rho(t,x)\right)\left(\rho(t,x+\eta)-\rho^*\right)f'\left(\eta^{-1}\int_x^{x+\eta}\rho(t,s)ds\right)g\left(\rho^*-\int_{x-1}^{x}(x-s)\rho(t,s)ds\right)dx\\
&\quad-\eta^{-1}\int_0^1\left(\rho^*-\rho(t,x)\right)^2\left|f'\left(\eta^{-1}\int_x^{x+\eta}\rho(t,s)ds\right)\right|g\left(\rho^*-\int_{x-1}^{x}(x-s)\rho(t,s)ds\right)dx\\
&\quad-\int_0^1\left(\rho(t,x)-\rho^*\right)^2 f\left(\eta^{-1}\int_x^{x+\eta}\rho(t,s)ds\right)g'\left(\rho^*-\int_{x-1}^{x}(x-s)\rho(t,s)ds\right)dx
\end{aligned} \tag{5.110}$$

Using the facts that $g(\rho)\geq 1$ for all $\rho\geq 0$ and $\left|\rho^*-\rho(t,x)\right|\left|\rho(t,x+\eta)-\rho^*\right|\leq\frac{1}{2}\left(\rho(t,x)-\rho^*\right)^2+\frac{1}{2}\left(\rho^*-\rho(t,x+\eta)\right)^2$, we obtain from (5.110) for all $t\geq 0$:

$$\begin{aligned}
\dot{V}(t)&\leq\frac{1}{2\eta}\int_0^1\left(\rho(t,x+\eta)-\rho^*\right)^2\left|f'\left(\eta^{-1}\int_x^{x+\eta}\rho(t,s)ds\right)\right|g\left(\rho^*-\int_{x-1}^{x}(x-s)\rho(t,s)ds\right)dx\\
&\quad-\frac{1}{2\eta}\int_0^1\left(\rho^*-\rho(t,x)\right)^2\left|f'\left(\eta^{-1}\int_x^{x+\eta}\rho(t,s)ds\right)\right|g\left(\rho^*-\int_{x-1}^{x}(x-s)\rho(t,s)ds\right)dx\\
&\quad-\int_0^1\left(\rho(t,x)-\rho^*\right)^2 f\left(\eta^{-1}\int_x^{x+\eta}\rho(t,s)ds\right)g'\left(\rho^*-\int_{x-1}^{x}(x-s)\rho(t,s)ds\right)dx
\end{aligned} \tag{5.111}$$

Using (2.16) and the fact that $\rho^*=\int_x^{x+1}\rho(t,s)ds$ for all $t\geq 0$, $x\in\Re$, we obtain the following estimates for all $t\geq 0$, $x\in\Re$:

$$\max\left(\rho^*-\frac{\rho_{\max}}{2},\frac{\rho_{\min}}{2}\right)\leq\rho^*-\int_{x-1}^{x}(x-s)\rho(t,s)ds\leq\min\left(\rho^*-\frac{\rho_{\min}}{2},\frac{\rho_{\max}}{2}\right) \tag{5.112}$$

$$\rho_{\min}\leq\eta^{-1}\int_x^{x+\eta}\rho(t,s)ds\leq\min\left(\eta^{-1}\rho^*,\rho_{\max}\right) \tag{5.113}$$

Consequently, using the fact that $\int_0^1 \left(\rho(t,x+\eta)-\rho^*\right)^2 dx = \int_0^1 \left(\rho(t,x+\eta)-\rho^*\right)^2 dx$ (a consequence of periodicity of $\rho[t]$), we obtain from (5.111), (5.112), (5.113), (3.3), (3.4), (3.5), (3.6), (3.7), (3.8) for all $t \geq 0$:

$$\dot{V}(t) \leq \left(\frac{1}{2\eta} F_{\max} g_{\max} - \frac{1}{2\eta} F_{\min} g_{\min} - f_{\min} G_{\min}\right) \int_0^1 \left(\rho(t,x)-\rho^*\right)^2 dx \tag{5.114}$$

Notice that the inequality $\frac{1}{2\rho_{\max}}\left(\rho-\rho^*\right)^2 \leq \rho \ln\left(\frac{\rho}{\rho^*}\right) + \rho^* - \rho \leq \frac{1}{2\rho_{\min}}\left(\rho-\rho^*\right)^2$ holds for all $\rho \in \left[\rho_{\min}, \rho_{\max}\right]$ with $0 < \rho_{\min} \leq \rho^* \leq \rho_{\max}$. Therefore, definition (5.109) implies the estimate:

$$\frac{1}{2\rho_{\max}} \int_0^1 \left(\rho(t,x)-\rho^*\right)^2 dx \leq V(t) \leq \frac{1}{2\rho_{\min}} \int_0^1 \left(\rho(t,x)-\rho^*\right)^2 dx \tag{5.115}$$

Combining (5.114) and (5.115), obtain the following differential inequality for all $t \geq 0$:

$$\dot{V}(t) \leq -\bar{c} V(t) \tag{5.116}$$

where $\bar{c} := \eta^{-1} \rho_{\min} \left(2\eta f_{\min} G_{\min} - F_{\max} g_{\max} + F_{\min} g_{\min}\right)$. The differential inequality (5.116) implies that $V(t) \leq \exp(-\bar{c}t) V(0)$ for all $t \geq 0$. The previous inequality in conjunction with (5.115) implies estimate (3.9). The proof is complete. ⊲

## 6. Concluding Remarks

The paper provided indications about the stabilizing effect of nudging in a ring-road when nudging is expressed by means of (2.3). However, Lemma 2.2 and the proof of Theorem 2.4 indicate that it is also possible to study more complicated speed adjustment feedback laws of the form

$$v(t,x) = \sum_{j=1}^{m} f_j\left(\int_x^{x+\eta_j} \omega_j(s-x)\rho(t,s)ds\right) g_j\left(\int_{x-\zeta_j}^{x} \tilde{\omega}_j(x-s)\rho(t,s)ds\right), \text{ for } t \geq 0, x \in \Re \tag{6.1}$$

where $\zeta_j, \eta_j > 0$ ($j = 1,...,m$) are constants, $g_j : \Re_+ \to \Re_+$ ($j = 1,...,m$) are non-decreasing, bounded functions, and $f_j : \Re_+ \to \Re_+$, $\omega_j : \Re_+ \to \Re_+$, $\tilde{\omega}_j : \Re_+ \to \Re_+$ ($j = 1,...,m$) are non-increasing functions. Such a research direction may allow the development of global stabilization results for ring-roads. This development has to be combined with the construction of appropriate Lyapunov functionals for the system.

Another research direction is the study of the effect of boundary conditions in non-local conservation laws on bounded domains. It is (in principle) possible to combine boundary feedback stabilization approaches with nudging and obtain even better results.

## Acknowledgments

The research leading to these results has received funding from the European Research Council under the European Union’s Horizon 2020 Research and Innovation programme/ ERC Grant Agreement n. [833915], project TrafficFluid.

# Appendix

**Proof of Lemma 2.2:** The proof for $(K+G)\in C^0(Per(\Re);Per(\Re))$ and $(\lambda K)\in C^0(Per(\Re);Per(\Re))$ is trivial. We focus on the proof for $(KG)\in C^0(Per(\Re);Per(\Re))$. Inequalities (2.8), (2.10), (2.13) can be shown easily. Inequality (2.9) for $K_N(\rho)G_N(\rho)$ is a direct consequence of the following equality

$$\begin{aligned}
&K_N(\rho)G_N(\rho)-K_N(\rho^{(1)})G_N(\rho^{(1)})-K_N(\tilde{\rho})G_N(\tilde{\rho})+K_N(\tilde{\rho}^{(1)})G_N(\tilde{\rho}^{(1)})\\
&=\left(K_N(\rho)-K_N(\rho^{(1)})-K_N(\tilde{\rho})+K_N(\tilde{\rho}^{(1)})\right)G_N(\rho)\\
&+\left(K_N(\rho^{(1)})-K_N(\tilde{\rho}^{(1)})\right)\left(G_N(\rho)-G_N(\rho^{(1)})\right)\\
&+\left(K_N(\tilde{\rho})-K_N(\tilde{\rho}^{(1)})\right)\left(G_N(\rho)-G_N(\tilde{\rho})\right)\\
&+K_N(\tilde{\rho}^{(1)})\left(G_N(\tilde{\rho}^{(1)})+G_N(\rho)-G_N(\rho^{(1)})-G_N(\tilde{\rho})\right)
\end{aligned}$$

inequalities (2.8), (2.9) (2.10) for $K_N(\rho),G_N(\rho)$ and the fact that $K_N(\rho),G_N(\rho)$ are bounded mappings. Similarly, inequality (2.11) for $K_N(\rho)G_N(\rho)$ is a consequence of the following equality

$$\begin{aligned}
&2K_N(\rho^{(1)})G_N(\rho^{(1)})-K_N(\rho)G_N(\rho)-K_N(\rho^{(2)})G_N(\rho^{(2)})\\
&=K_N(\rho^{(1)})\left(2G_N(\rho^{(1)})-G_N(\rho)-G_N(\rho^{(2)})\right)\\
&+G_N(\rho^{(2)})\left(2K_N(\rho^{(1)})-K_N(\rho^{(2)})-K_N(\rho)\right)\\
&+\left(K_N(\rho)-K_N(\rho^{(1)})\right)\left(G_N(\rho^{(2)})-G_N(\rho)\right)
\end{aligned}$$

as well as inequalities (2.10), (2.11) for $K_N(\rho),G_N(\rho)$ and the fact that $K_N(\rho),G_N(\rho)$ are bounded mappings. Finally, inequality (2.12) for $K_N(\rho)G_N(\rho)$ is a consequence of the following equality

$$\begin{aligned}
&3K_N(\rho^{(1)})G_N(\rho^{(1)})+K_N(\rho^{(3)})G_N(\rho^{(3)})-K_N(\rho)G_N(\rho)-3K_N(\rho^{(2)})G_N(\rho^{(2)})\\
&=\left(3K_N(\rho^{(1)})+K_N(\rho^{(3)})-K_N(\rho)-3K_N(\rho^{(2)})\right)G_N(\rho^{(1)})\\
&+\left(K_N(\rho^{(3)})-2K_N(\rho^{(2)})+K_N(\rho^{(1)})\right)\left(G_N(\rho^{(3)})-G_N(\rho^{(1)})\right)\\
&+\left(K_N(\rho^{(2)})-K_N(\rho^{(1)})\right)\left(G_N(\rho^{(3)})-2G_N(\rho^{(2)})+G_N(\rho^{(1)})\right)\\
&+2\left(K_N(\rho^{(2)})-K_N(\rho^{(1)})\right)\left(G_N(\rho^{(2)})-2G_N(\rho^{(1)})+G_N(\rho)\right)\\
&+\left(K_N(\rho^{(2)})-2K_N(\rho^{(1)})+K_N(\rho)\right)\left(G_N(\rho^{(1)})-G_N(\rho)\right)\\
&+K_N(\rho^{(2)})\left(3G_N(\rho^{(1)})+G_N(\rho^{(3)})-G_N(\rho)-3G_N(\rho^{(2)})\right)
\end{aligned}$$

as well as inequalities (2.10), (2.11), (2.12) for $K_N(\rho), G_N(\rho)$ and the fact that $K_N(\rho), G_N(\rho)$ are bounded mappings. The proof is complete. $\triangleleft$